\documentclass{article}
\usepackage{amsmath,amsfonts,amssymb}
\usepackage{xy}

\title{Fine-Wilf graphs and the generalized Fine-Wilf theorem}

\author{Stuart A. Rankin}

\begin{document}
\maketitle
\begin{abstract}
 In 1962, R. C. Lyndon and 
 M. P. Shutzenberger established that for any positive integers $r$ and $s$, any sequence of length at least $r+s$ 
 that is both $r$-periodic and $s$-periodic is then $(r,s)$-periodic. Shortly thereafter (1965), N. J. Fine and H. S. Wilf
 proved that for any positive integers $r$ and $s$, if $a$ is an infinite seqeunce of period $r$
 and $b$ is an infinite sequence of period $s$ such that $a_i=b_i$ for all $i$ with $1\le i\le r+s-(r,s)$,
 then $a=b$. This is equivalent to the following result, which is commonly referred to as the Fine-Wilf
 theorem: for any positive integers $r$ and $s$, if $w$ is a finite sequence that is both $r$-periodic and $s$-periodic,
 and $|w|\ge r+s-(r,s)$, then $w$ is $(r,s)$-periodic. Fine and Wilf also asserted that this bound is best
 possible, in the sense that for any positive integers $r$ and $s$, there exists a word $w$ of length $r+s-(r,s)-1$ that
 is both $r$-periodic and $s$-periodic, but not $(r,s)$-periodic. This sharpness result has since been established,
 and these extremal sequences are now much studied. Among other results, it is known that for a given $r$ and $s$, there
 is a unique (up to relabelling) sequence of length $r+s-(r,s)-1$ that is both $r$-periodic and $s$-periodic, but
 not $(r,s)$-periodic, and in this sequence, exactly two distinct entries appear. 
 
 The Fine-Wilf theorem was generalized to finite sequences 
 with three periods by M. G. Castelli, F. Mignosi, and A. Restivo. They introduced a function $f$ from the
 set of all ordered triples of nonnegative integers to the set of positive integers with the property that if 
 $w$ is a finite sequence with periods $p_1$, $p_2$, and $p_3$, and $|w|\ge f(p)$, where $p=(p_1,p_2,p_3)$, then 
 $w$ is $(p)$-periodic as well. They further established a condition on $p$ under which the bound $f(p)$ is best 
 possible. In support of their work, they introduced the graphs that we shall refer to as Fine-Wilf graphs.
 The work of Castelli et al. was generalized by J. Justin, and more broadly by R. Tijdeman and L. Zamboni, 
 who introduced a function $fw$ from the set of all sequences of nonnegative integers to the set of positive 
 integers, and they proved that for a sequence $p=(p_1,p_2,\ldots,p_n)$, a finite sequence
 $w$ with periods $p_i$, $i=1,2,\ldots, n$ and length at least $fw(p)$ must be $(p)$-periodic as well, 
 and that there exists a sequence $w$ of length $fw(p)-1$ that is $p_i$-periodic for all $i$, but not
 $(p)$-periodic.
 
 In this paper, we follow ideas introduced by S. Constantinescu and L. Ilie to obtain an alternative
 formulation of $f$ and $fw$, and we use this formulation to establish important properties of $f$
 and $fw$, obtaining in particular new upper and lower bounds for each. We also begin an investigation 
 of Fine-Wilf graphs for arbitrary finite sequences with a view to understanding how the graph
 may be used to better understand $f$ and $fw$.
\end{abstract}

\def\z{\mathbb{Z}}
\def\zplus{\z^+}
\def\ofs{\textrm{OFS}(\zplus)}
\def\gcdp#1{\textrm{gcd}(#1)}
\def\np{\hbox{$+$}}
\def\nm{\hbox{$-$}}
\def\aut{\text{Aut}}
\def\seq#1{\{#1\}}
\def\set{\{\mkern.8mu }
\def\endset{\mkern.8mu\}} 
\def\rest#1{\,\hbox{\vrule height 6pt depth 5pt width .5pt}_{\,#1}}            
\def\from{\mkern2mu\hbox{\rm :}\mkern2mu}
\def\nat{\mathbb{N}}
\def\component#1{\kappa_{#1}}
\def\table#1\endtable{\begin{tabular}#1 \end{tabular}}
\newtheorem{theorem}{Theorem}
\newtheorem{definition}[theorem]{Definition}
\newtheorem{lemma}[theorem]{Lemma}
\newtheorem{corollary}[theorem]{Corollary}
\newtheorem{proposition}[theorem]{Proposition}
\newtheorem{example}[theorem]{Example}
\def\proof{\ifdim\lastskip<\baselineskip\relax\removelastskip
  \vskip\baselineskip\fi\leavevmode
  {\it Proof. }}
\def\endproof{\par\medskip}
\def\strutdepth{\dp\strutbox}
\def\endproofmarker{\strut\vadjust{\kern-2\strutdepth\epmarker}}
\def\epmarker{\vbox to \strutdepth{\baselineskip\strutdepth\vss\hfill{%
\hbox to 0pt{\hss\vrule height 4pt width 4pt depth 0pt}\null}}}
\let\oldendproof=\endproof
\def\endproof{\endproofmarker\oldendproof}
\let \cong\equiv


\section{Introduction} 
 For any positive integer $r$, a finite sequence $w=(a_1,a_2,\ldots, a_n)$ is said to have period $r$, or to be $r$-periodic,
 if for every positive integer $i$ for which $i,i+r\le n$, $a_i=a_{i+r}$. In 1962, R.\ C.\ Lyndon and 
 M.\ P.\ Shutzenberger \cite{lyndon} established that for any positive integers $r$ and $s$, if $w$ is both $r$-periodic and 
 $s$-periodic, and $|w|\ge r+s$, then $w$ is $\gcdp{r,s}$-periodic. Shortly thereafter (1965), N. J. Fine and H. S. Wilf
 \cite{fine} proved that for any positive integers $r$ and $s$, if $\seq{a_i}$ is an infinite seqeunce of period $r$
 and $\seq{b_i}$ is an infinite sequence of period $s$ such that $a_i=b_i$ for all $i$ with $1\le i\le r+s-\gcdp{r,s}$,
 then $a_i=b_i$ for all $i$. This is equivalent to the following result, which is commonly referred to as the Fine-Wilf
 theorem: for any positive integers $r$ and $s$, if $w$ is a finite sequence that is both $r$-periodic and $s$-periodic,
 and $|w|\ge r+s-\gcdp{r,s}$, then $w$ is $\gcdp{r,s}$-periodic. It was also asserted in \cite{fine} that this bound is best
 possible, in the sense that for any positive integers $r$ and $s$, there exists a word $w$ of length $r+s-\gcdp{r,s}-1$ that
 is both $r$-periodic and $s$-periodic, but not $\gcdp{r,s}$-periodic. This sharpness result has since been established,
 and these extremal sequences are now much studied. Among other results, it is known that for a given $r$ and $s$, there
 is a unique (up to relabelling) sequence of length $r+s-\gcdp{r,s}-1$ that is both $r$-periodic and $s$-periodic, but
 not $\gcdp{r,s}$-periodic, and in this sequence, exactly two distinct entries appear. For example, for $r=2$ and $s=3$,
 the sequence is $(0,1,0)$. 
 
 Nearly thirty-five years later (1999), the Fine-Wilf theorem was generalized to finite sequences 
 with three periods by M. G. Castelli, F. Mignosi, and A. Restivo \cite{castelli}. They introduced a function $f$ from the
 set of all ordered triples of nonnegative integers to the set of positive integers with the property that if $w$ is a finite sequence
 with periods $p_1$, $p_2$, and $p_3$, and $|w|\ge f(p)$, where $p=(p_1,p_2,p_3)$, then $w$ is $\gcdp{p}$-periodic as well. They further
 established a condition on $p$ under which the bound $f(p)$ is best possible. The sequences $p$ that
 met this condition were precisely those for which the unique (up to relabelling) finite sequence of greatest length and with the
 greatest possible number of distinct entries that had
 periodicity $p_1$, $p_2$, and $p_3$, but not $\gcdp{p_1,p_2,p_3}$ had exactly three distinct entries.  In support of their work,
 they introduced
 the graphs that we shall refer to as Fine-Wilf graphs $G(p_1,p_2,p_3,n)$, where $p_1$, $p_2$, and $p_3$ are distinct 
 nonnegative integers, $n$ is a positive integer, and $G(p_1,p_2,p_3,n)$ denotes the graph with vertex set 
 $\set 1,2,\ldots,n\endset$, and edge set
 $$
   \set \set i,j\endset\mid |i-j|\in \set p_1,p_2,p_3\endset\endset.
 $$   

 The work of Castelli et al. was followed immediately (2000) by work of
 J. Justin \cite{justin}, who extended the definition of the function
 $f$ to all finite sequences of nonnegative integers, with analagous
 results.
 
 A broader generalization of the work of Castelli et al. was then given
 by R. Tijdeman and L. Zamboni \cite{tij} (2003). They introduced a
 function,  which we shall denote as $fw$, from the set of all sequences
 of nonnegative integers to the set of positive integers, and they
 proved that for a sequence $p=(p_1,p_2,\ldots,p_n)$, a finite sequence
 $w$ with periods $p_i$, $i=1,2,\ldots, n$ and length at least $fw(p)$
 must be $\gcdp{p}$-periodic as well, and that there exists a sequence
 $w$ of length $fw(p)-1$ that is $p_i$-periodic for all $i$, but not
 $\gcdp{p}$-periodic. At nearly the same time (2005), and independently
 of the work of Tijdeman and Zamboni, S. Constantinescu and L. Ilie
 \cite{ilie}  described what amounts to an extension of the function $f$
 of Castelli et al., and used $f$ to compute a related function that
 gives the best  bound in all cases. Of course, this related function is
 the function $fw$, but the evalution of $fw$ as described by
 Constantinescu and Ilie is quite different from that of Tijdeman and
 Zamboni.
 
 In this paper, we establish important properties of the functions $f$
 and $fw$. In particular, we introduce new upper and lower  bounds for
 $f$. We also begin an investigation of Fine-Wilf graphs for arbitrary
 $p_1,p_2,\ldots,p_n$, with a view to understanding how the graph
 depends on the values $p_1,p_2,\ldots,p_n$. 

\section{Generalization of the Fine-Wilf theorem}
Let $\ofs$ denote the set of all strictly increasing finite sequences of positive integers. For
$p\in \ofs$, let $\gcdp{p}$ denote the greatest common divisor of the entries in $p$, let $|p\mkern1mu|$ 
denote the length of $p$, and for $1\le i\le |p\mkern1mu|$, let $p_i$ denote the $i^{th}$ entry of $p$ and
let $p\rest i$ denote the truncated sequence $(p_1,p_2,\ldots,p_i)$. Finally, let $\max(p)$ and $\min(p)$
denote $p_n$, respectively $p_1$, where $n=|p\mkern1mu|$.

\begin{definition}
 Let $R\from \ofs\to \ofs$ denote the function defined as follows. For $p\in \ofs$, 
 $R(p)= p$ if $|p\mkern1mu|=1$. If $n=|p\mkern1mu|>1$, then form the sequence $(p_2-p_1,p_3-p_1,\ldots,p_n-p_1)\in \ofs$,
 and, if $p_1$ does not appear in this sequence, insert $p_1$ so the result is an element of $\ofs$. The
 sequence that results (of length either $n-1$ or $n$) is denoted by $R(p)$. Moreover, we shall define
 $p^{(i)}\in\ofs$ for $i\ge 0$ as follows: $p^{(0)}=p$, and for $k\ge0$, $p^{(k+1)}=R(p^{(k)})$. 
\end{definition}

Note that for any $p\in\ofs$, $\gcdp{p}=\gcdp{R(p)}$.

\begin{definition}
 Define $f\from\ofs\to\zplus$ by induction on $\max(p)$ as follows. If $p\in\ofs$ has 
 $\max(p)=1$, then $p=(1)$, and we define $f((1))=1$. Then for $p\in \ofs$ with $\max(p)>1$,
 define
 $$
 f(p)=\begin{cases} p_1 &\text{if $|p\mkern1mu|=1$}\\
             p_1+f(R(p)) &\text{if $|p\mkern1mu|>1$.}\end{cases}
 $$	     
 Moreover, the column of sequences whose $i^{th}$ row is $p^{(i)}$, $i\ge 0$, and whose last row is $p^{(m)}$,
 where $m$ is least subject to the requirement that $|p^{(m)}|=1$ shall be called the tableau for the calculation
 of $f(p)$.
\end{definition}

For example, if $p=(4,7)$ or $(4,7,9)$, then $f(p)=10$, as can be seen from the tableaux below.

\vskip-\baselineskip
{\centering\table{c@{\hskip20pt}c}
\table{l}
4,7\\
3,4\\
1,3\\
1,2\\
1\endtable
&
\table{l}
4,7,9\\
3,4,5\\
1,2,3\\
1,2\\
1\endtable\\
\noalign{\vskip0pt}
\multicolumn{1}{c}{\table{c}Tableau for the\\
calculation of $f(p)$ for\\
$p=(4,7)$.\endtable}
&
\table{c}Tableau for the\\
calculation of $f(p)$ for\\
$p=(4,7,9)$.\endtable
\endtable
\par}

\begin{lemma}\label{f bounds}
 For any $p\in\ofs$, $f(p)\ge \max(p)$. Furthermore, if $|p\mkern1mu|>1$, then $f(p)\ge 2p_1$.
\end{lemma}

\proof
 By induction on $\max(p)$. It is certainly true when $p=(1)$, and this is the base case $\max(p)=1$.
 Suppose now that $m>1$ is an integer such that the result holds for $p\in\ofs$ with $1\le \max(p)<m$, and
 let $p\in\ofs$ be such that $\max(p)=m$. If $|p\mkern1mu|=1$, then $f(p)=p_1=\max(p)$, so the result holds trivially.
 Suppose that $|p\mkern1mu|>1$. Then $f(p)=p_1+f(R(p))$. By hypothesis, $f(R(p))\ge \max(R(p))$. If
 $\max(p)-p_1>p_1$, then $\max(R(p))=\max(p)-p_1$, otherwise $\max(R(p))=p_1$. In the
 former case, we have $f(p)=p_1+f(R(p))\ge p_1+\max(p)-p_1=\max(p)>2p_1$, while in the latter case, we
 have $f(p)=p_1+f(R(p))\ge p_1+p_1=2p_1$, and this case occurs when $\max(p)-p_1\le p_1$, or $2p_1\ge \max(p)$.
 Thus in either case, we have $f(p)\ge \max(p)$, and (recall that $|p\mkern1mu|>1$ in these cases), we have
 $f(p)\ge 2p_1$. The result follows now by induction.
\endproof

\begin{corollary}\label{frp}
 For $p\in\ofs$ with $|p\mkern1mu|>1$, $f(R(p))\ge p_1$.
\end{corollary}

\proof
 By Lemma \ref{f bounds}, $f(p)=p_1+f(R(p))\ge 2p_1$, so $f(R(p))\ge p_1$.
\endproof

\begin{definition}
 Define $fw\from\ofs\to\zplus$ as follows. Let $p\in\ofs$. Then
 $fw(p)=fw(p\rest{|p\mkern1mu|-1})$ if $|p\mkern1mu|>1$, $\gcdp{p\rest {|p\mkern1mu|-1}}=\gcdp{p}$, and 
 $\max(p)\ge f(p\rest{|p\mkern1mu|-1})$, otherwise $fw(p)=f(p)$.
\end{definition}

We shall show later (see Proposition \ref{truncated sequence fw}) that if $p\in\ofs$ with $|p\mkern1mu|>1$, 
$\gcdp{p\rest {|p\mkern1mu|-1}}=\gcdp{p}$, and 
$\max(p)<fw(p\rest{|p\mkern1mu|-1})$, then $fw(p)\le fw(p\rest{|p\mkern1mu|-1})$.

\begin{proposition}\label{fw lower bound}
 Let $p\in\ofs$ with $p_1\ne \gcdp{p}$. Then $fw(p)\ge 2p_1$.
\end{proposition}

\proof
 Since $\gcdp{p}\ne p_1$, $|p\mkern1mu|\ge 2$. If $fw(p)=f(p)$, then the result follows from
 Lemma \ref{f bounds}. Suppose that $fw(p)\ne f(p)$. Then (since $\gcdp{p}\ne p_1$),
 there exists an index $i>1$ such that $fw(p)=f(p\rest{i})$, and by Lemma \ref{f bounds},
 $f(p\rest{i})\ge 2p_1$.
\endproof
 
\begin{proposition}\label{reduced}
 Let $p\in \ofs$, $d=\gcdp{p}$, and $\displaystyle\frac{p}{d}= (\frac{p_1}{d},\ldots,\frac{p_{|p\mkern1mu|}}{d})$. 
 Then $\displaystyle\frac{p}{d}\in \ofs$, $\displaystyle f(p)=d\, f(\frac{p}{d})$, 
 and $\displaystyle fw(p)=d\,fw(\frac{p}{d})$.
\end{proposition}

\proof
 By induction on $\max(p)$. For $p\in\ofs$ such that $\max(p)=1$, it must be that $p=(1)$
 and $d=1$, so the result holds in this case. Suppose now that $p\in\ofs$ with $\max(p)>1$ and
 that the result holds for all sequences in $\ofs$ with smaller maximum entry.
 If $|p\mkern1mu|=1$, then $d=p_1$, $p=(p_1)$, and $p/d=(1)$, so $d\,f(p/d)=p_1=f((p_1))=f(p)$ and 
 $d\,fw(p/d)=p_1=fw((p_1))=fw(p)$. Suppose that $n=|p\mkern.8mu|>1$, so that $n=|p/d\mkern.8mu|$ as well. 
 Since there is nothing to prove if $d=1$, suppose that $d>1$.
 We have  $f(p)=p_1+f(R(p))$, and since $\gcdp{R(p)}=
 \gcdp{p}$, it follows from our induction hypothesis that $f(R(p))=d\,f(R(p)/d)=d\,f(R(p/d))$ and thus 
 $f(p)=p_1+d\,f(R(p/d))=d(p_1/d+f(R(p/d)))=d\,f(p/d)$. As for $fw$, $fw(p)=f(p)$ if
 $\gcdp{p}\ne \gcdp{p\rest{n-1}}$, or $\gcdp{p}= \gcdp{p\rest{n-1}}$ and
 $p_n\ge f(p\rest{n-1})$, otherwise $fw(p)=fw(p\rest{n-1})$.  
 Observe that $\gcdp{p}= \gcdp{p\rest{n-1}}$ if and only if 
 $\gcdp{p/d}=1= \gcdp{p\rest{n-1}/d}=\gcdp{(p/d)\rest{n-1}}$.
 Suppose first that $\gcdp{p}\ne \gcdp{p\rest{n-1}}$, in which case we have
 $\gcdp{p/d}=1\ne \gcdp{p\rest{n-1}/d}=\gcdp{(p/d)\rest{n-1}}$ and thus 
 $fw(p)=f(p)$ and $fw(p/d)=f(p/d)$, so $fw(p)=f(p)=d\,f(p/d)=d\,fw(p/d)$. Now suppose that
 $\gcdp{p}= \gcdp{p\rest{n-1}}$, so  $\gcdp{p/d}=\gcdp{(p/d)\rest{n-1}}$.
 If $p_n<f(p\rest{n-1})$, then $fw(p)=f(p)=d\,f(p/d)$ and $p_n/d<f((p/d)\rest{n-1})$, so $d\,fw(p/d)=d\,f(p/d)=f(p)=fw(p)$
 in this case. Finally, suppose that $p_n\ge f(p\rest{n-1})$, so $fw(p)=fw(p\rest{n-1})$, and
 $p_n/d\ge f((p/d)\rest{n-1})$. Thus $fw(p/d)=fw((p/d)\rest{n-1})$. Since $\gcdp{p\rest{n-1}}=d>1$, $\max(p/d)=p_n/d<p_n=\max(p)$,
 and so we may apply the induction hypothesis to $(p/d)\rest{n-1}$ to obtain $fw(p\rest{n-1})=d\,fw((p/d)\rest{n-1})$.
 Thus $fw(p)=fw(p\rest{n-1})=d\,fw((p/d)\rest{n-1})=d\,fw(p/d)$, as required.
 Since in each case we have $fw(p)=d\,fw(p/d)$, the result follows by induction.
\endproof

The next result gives an important lower bound for $f(p)$, and this result can be viewed in a sense
as a generalization of the Fine-Wilf theorem. We will later obtain an upper bound (see Proposition \ref{upper bound for f},
also Proposition \ref{upper bound for fw(p)}) for $f(p)$ and the
combination of that upper bound with the following lower bound, when applied in the case of $|p\mkern1mu|=2$,
will give the Fine-Wilf theorem.

\begin{proposition}\label{good lower bound for f}
 Let $p\in\ofs$ with $|p\mkern1mu|>1$. Then 
 $$
   f(p)\ge \frac{-\gcdp{p}+\sum_{i=1}^{|p\mkern1mu|} p_i}{|p\mkern1mu|-1},
 $$
 and if equality holds, but for some $i$, $p_i=2p_1$, then $i=|p\mkern1mu|$ and $f(p)=2p_1=\max(p)$.
\end{proposition}

\proof
 By Proposition \ref{reduced}, it suffices to prove the result only for $p$ with $\gcdp{p}=1$. We prove by induction on $\max(p)$
 that $|p\mkern1mu|>1$ and $\gcdp{p}=1$ implies $\displaystyle f(p)\ge \frac{-1+\sum_{i=1}^{|p\mkern1mu|} p_i}{|p\mkern1mu|-1}$. The base case,
 $\max(p)=1$, is trivially true, so suppose that $p\in\ofs$ has $\max(p)>1$ and that the result holds for all elements of
 $\ofs$ with smaller maximum entry. Further suppose that $|p\mkern1mu|>1$ and $\gcdp{p}=1$, and let $n=|p\mkern1mu|$. Consider first the case when
 $p_i\ne 2p_1$ for every $i=1,2,\ldots,n$. Then $|R(p)|=|p\mkern1mu|=n>1$, $\gcdp{R(p)}=\gcdp{p}=1$, and the $n$ entries of
 $R(p)$ are $p_1,p_2-p_1,p_3-p_1,\ldots,p_n-p_1$ (with $p_1$ not necessarily in the correct position). We may apply
 the induction hypothesis to obtain that
 $$
  f(p)=p_1+f(R(p))\ge p_1+\frac{-1+p_1+\sum_{i=2}^n(p_i-p_1)}{n-1}=\frac{-1+\sum_{i=1}^{|p\mkern1mu|} p_i}{|p\mkern1mu|-1},
 $$
 as required.
 Suppose now that for some $j$, $p_j=2p_1$. Then $|R(p)|=|p\mkern1mu|-1=n-1$ and $R(p)=
 (p_2-p_1,\ldots,p_{j-1}-p_1,p_1,p_{j+1}-p_1,\ldots,p_n-p_1)$. If $n=2$, then since $\gcdp{p}=1$, $p=(1,2)$ and $f(p)=2$,
 while $(1+2-1)/1=2$, so the result holds in this case. Note that in this case we have equality, and $f(p)=2p_1=\max(p)$,
 as required. Otherwise, $n>2$, so $|R(p)|=n-1>1$, $\gcdp{R(p)}=\gcdp{p}=1$, and
 we may apply the induction hypothesis to obtain that 
 $$
  f(p)\ge p_1+\frac{-1+p_1+\sum_{i=2,\ i\ne j}^n(p_i-p_1)}{n-2}=\frac{-1-2p_1+\sum_{i=1}^np_i}{n-2}.
 $$
 Now $\displaystyle\frac{-1-2p_1+\sum_{i=1}^np_i}{n-2}\ge \frac{-1+\sum_{i=1}^n p_i}{n-1}$ if and only if
 $$
  -1+\sum_{i=2,i\ne j}^n p_i\ge (2n-5)p_1.
 $$
 Suppose that $-1+\sum_{i=2,i\ne j}^n p_i\ge (2n-5)p_1$. Then $f(p)\ge \frac{-1-2p_1+\sum_{i=1}^np_i}{n-2}\ge \frac{-1+\sum_{i=1}^n p_i}{n-1}$.
 Further suppose that $f(p)= \frac{-1+\sum_{i=1}^n p_i}{n-1}$. Then $\frac{-1-2p_1+\sum_{i=1}^np_i}{n-2}= \frac{-1+\sum_{i=1}^n p_i}{n-1}$
 and so $\frac{-1+\sum_{i=1}^n p_i}{n-1}=2p_1$. But then $2p_1=f(p)\ge p_n\ge p_j=2p_1$, so $p_j=p_n$ and thus $j=n$, as
 required.
 
 Now suppose that $-1+\sum_{i=2,i\ne j}^n p_i < (2n-5)p_1$. Then we have $-1+\sum_{i=1}^n p_i=-1+3p_1+\sum_{i=2,i\ne j}^n p_i\le (2n-2)p_1$, and so
 $\frac{-1+\sum_{i=1}^n p_i}{n-1}\le 2p_1$. By Lemma \ref{f bounds}, $f(p)\ge 2p_1$, so we have
 $\displaystyle f(p)\ge 2p_1\ge \frac{-1+\sum_{i=1}^n p_i}{n-1}$. Moreover, if $\displaystyle f(p)=\frac{-1+\sum_{i=1}^n p_i}{n-1}$,
 then $f(p)=2p_1$ and as before, $2p_1=f(p)\ge p_n\ge p_j=2p_1$, and so $j=n$. This completes the proof of the inductive
 step and so the result follows.
\endproof
 
\begin{proposition}\label{min(p)=gcd(p)}
 Let $p\in\ofs$ with $\min(p)=\gcdp{p}$. Then $fw(p)=\min(p)$ and $f(p)=\max(p)$.
\end{proposition}

\proof
 By induction on $|p\mkern1mu|$. If $|p\mkern1mu|=1$, then $fw(p)=f(p)=p_1=\min(p)=\max(p)$. Suppose now
 that $n=|p\mkern1mu|>1$. Since $\gcdp{p}=p_1=\gcdp{p\rest{n-1}}$ and by hypothesis, $f(p\rest{n-1})=
 p_{n-1}<p_n$, we have $fw(p)=fw(p\rest{n-1})=p_1$, while $f(p)=p_1+f(R(p))$. Since every term
 is a multiple of $p_1$, $\max(R(p))=\max(p)-p_1$, so by hypothesis, we have $f(p)=p_1+\max(p)-p_1=\max(p)$.
\endproof
 
\begin{proposition}\label{ready for optimality}
 For any $p\in \ofs$, the following hold.
 \begin{enumerate}
  \item  $f(p)\ge fw(p)$. 
  \item If $|p\mkern1mu|\ge 2$, and $i$ is such that $1\le i<|p\mkern1mu|$ and $\gcdp{p}=\gcdp{p\rest{j}}$ for all $i\le j< |p\mkern1mu|$, and 
 $p_{i+1}\ge f(p\rest{i})$, then $f(p\rest{j})=p_j$ for all $j$ with $i+1\le j\le |p\mkern1mu|$.
 \end{enumerate}
\end{proposition}

\proof
 By induction on $\max(p)$. Let $p\in \ofs$. If $\max(p)=1$, then $p=(1)$ and so $fw(p)=f(p)=1$.
 Suppose now that $p\in\ofs$ has $\max(p)>1$ and the result holds for all elements of $\ofs$
 with smaller maximum entry. If $|p\mkern1mu|=1$,
 then $fw(p)=f(p)$, so we may suppose that $n=|p\mkern1mu|\ge2$. If $\gcdp{p}\ne \gcdp{p\rest{n-1}}$
 or $\gcdp{p}= \gcdp{p\rest{n-1}}$ but $p_n< f(p\rest{n-1})$, then $fw(p)=f(p)$ by definition,
 so we may suppose that $\gcdp{p}=\gcdp{p\rest{n-1}}$ and $p_n\ge f(p\rest{n-1})$. 
 Let $i$ be such that $\gcdp{p}=\gcdp{p\rest{j}}$ for
 all $i\le j< n=|p\mkern1mu|$, and $p_{i+1}\ge f(p\rest{i})$. Furthermore, we may assume without loss of
 generality that $i$ is minimal in this regard. If $i=1$, then $p_1=\gcdp{p}=fw(p)$, and thus the
 result follows from Proposition \ref{min(p)=gcd(p)}. Consider now the case when $i>1$. Suppose first
 that $i<n-1$. Then by the inductive hypothesis applied to $p\rest{n-1}$, $f(p\rest{j})=p_j$ for 
 all $i+1\le j<n$. If $f(p)=p_n$, then, since by our induction hypothesis,
 $f(p\rest{n-1})=p_{n-1}<p_n$, we will have by definition that $fw(p)=fw(p\rest{n-1})$
 and so $f(p)=p_n>f(p\rest{n-1})\ge fw(p\rest{n-1})=fw(p)$, as required. It therefore suffices to prove
 that $f(p)=p_n$. In this case, $n>i>1$, $n-1>1$, and so we may apply Lemma \ref{f bounds} to 
 conclude that $f(p)\ge 2p_1$ and thus 
 $p_n> f(p\rest{n-1})\ge2p_1$. It follows that $p_n-p_1>p_1$. Now, $f(p)=p_1+f(R(p))$ and $p_n-p_1>p_1$,
 so $\max(R(p))=p_n-p_1$. Furthermore, we have $p_n-p_1>p_{n-1}-p_1=f(p\rest{n-1})-p_1=f(R(p\rest{n-1}))$.
 Since $\max(R(p))=p_n-p_1>p_1$, we have $R(p)\rest{|R(p)|-1}=R(p\rest{n-1})$, so 
 $$
   \max(R(p))=p_n-p_1>f(R(p\rest{n-1}))=f(R(p)\rest{|R(p)|-1}).
 $$
 Since $\max(R(p))<\max(p)$, our induction hypothesis applies to $R(p)$ and since 
 $$
   |R(p)|\ge2,\quad \gcdp{R(p)}=\gcdp{p}= \gcdp{p\rest{n-1}}=\gcdp{R(p\rest{n-1})}
 $$ 
 and $\max(R(p))>f(R(p)\rest{|R(p)|-1})$,
 we conclude that $f(R(p))=\max(R(p))=p_n-p_1$. Thus $f(p)=p_1+f(R(p))=p_1+p_n-p_1=p_n$.
 
 It remains to consider the case when $i=n-1$. We have $p_n\ge f(p\rest{n-1})= p_1+f(R(p\rest{n-1})$. Furthermore, we have
 $\gcdp{p}=\gcdp{p\rest{n-1}}$. We wish to 
 prove that $f(p)=p_n$ (and thus $f(p)\ge fw(p\rest{n-1})
 =fw(p)$ as well). Since $n>1$, we have by Lemma \ref{f bounds} that
 $p_n\ge f(p\rest{n-1})\ge 2p_1$, and so $p_n-p_1\ge p_1$.
 Suppose first that $p_n-p_1>p_1$. Then as above,
 we have $R(p)=(p_2-p_1,p_3-p_1,\ldots,p_1,\ldots,p_n-p_1)$, and so
 $R(p)\rest{|R(p)|-1}=R(p\rest{n-1})$. It then follows that
 $$
   \max(R(p))=p_n-p_1>f(R(p\rest{n-1})=f(R(p)\rest{|R(p)|-1}).
 $$
Then by our induction hypothesis (since $\gcdp{R(p)}=\gcdp{p}=\gcdp{p\rest{n-1}}=\gcdp{R(p\rest{n-1})}$),
 we have $f(R(p))=\max(R(p))=p_n-p_1$, which yields $f(p)=p_1+f(R(p))=p_1+p_n-p_1=p_n$.
 Finally, consider the case when $p_n-p_1=p_1$, or $p_n=2p_1$. In this case, we have
 $R(p)=(p_2-p_1,p_3-p_1,\ldots,p_{n-1}-p_1,p_1)=R(p\rest{n-1})$. As well, from
 $2p_1=p_n\ge f(p\rest{n-1})\ge 2p_1$, we obtain that $f(p\rest{n-1})=2p_1$, and
 thus $p_1+f(R(p\rest{n-1}))=2p_1$, or $f(R(p\rest{n-1}))=p_1$.
 Now, $f(p)=p_1+f(R(p))=p_1+f(R(p\rest{n-1}))=p_1+p_1=2p_1=p_n$, as required.
\endproof

In particular, if $p\in\ofs$ has $n=|p\mkern1mu|>1$, $\gcdp{p\rest{n-1}}=\gcdp{p}$, and $p_{n}\ge f(p\rest{n-1})$, then 
$f(p)=p_n$.

\begin{definition}\label{trim}
 Let $p\in\ofs$. We say that $p$ is trim if either $|p\mkern1mu|=1$, or else $n=|p\mkern1mu|>1$ and either
 $\gcdp{p}\ne \gcdp{p\rest{n-1}}$ or else $\gcdp{p}= \gcdp{p\rest{n-1}}$ but $\max(p)<f(p\rest{n-1})$. 
 For any $p\in\ofs$, there exists $i$ with $1\le i\le |p|$ such
 that $p\rest{i}$ is trim, and $q\in\ofs$ is called the trimmed form of $p$ if $q=p\rest{i}$ where $i$ is maximal
 with respect to the property $p\rest{i}$ is trim. 
\end{definition}
 
We note that even if $p$ is trim, there may exist $i$ with $1<i<|p|$ such that $p\rest{i}$ is not trim.

\begin{corollary}\label{not trim}
 Let $p\in\ofs$. If $p$ is not trim, then $f(p)=\max(p)$.
\end{corollary}

\proof
 Since $|p\mkern.8mu|=1$ implies that $p$ is trim, we have $n=|p\mkern.8mu|>1$. Since $\gcdp{p}=\gcdp{p\rest{n-1}}$
 and $p_n\ge f(p\rest{n-1})$, we may take $i=n-1$ in Proposition \ref{ready for optimality} (ii) to obtain
 that $f(p)=p_n=\max(p)$.
\endproof 

For $p\in\ofs$ with $n=|p\mkern1mu|\ge 2$, if $\gcdp{p}= \gcdp{p\rest{n-1}}$ and $\max(p)\ge fw(p\rest{n-1})$, then we shall say
that $p\rest{n-1}$ is obtained by trimming $p$. Evidently, for any $p\in\ofs$, we may iteratively apply the trimming
operation to obtain $q=p\rest{j}$ for some $j>1$ with $q$ trim, and we note that $fw(p)=fw(q)$. 

\begin{lemma}\label{benefits of trim}
 Let $p\in\ofs$. If $p$ is trim, then $fw(p)=f(p)$. Furthermore, if $|p\mkern1mu|>1$, then $\min(p)>\gcdp{p}$.
\end{lemma}

\proof
 That $fw(p)=f(p)$ is immediate from the definition of $fw$. If $|p\mkern1mu|>1$, then by Proposition \ref{min(p)=gcd(p)}, 
 we have $\min(p)>\gcdp{p}$.
\endproof 

\begin{proposition}\label{lower bound for fw}
 Let $p\in\ofs$. If $p$ is trim with $|p\mkern1mu|>1$, then $f(p)> \max(p)$. In addition,
 if $R(p)$ is not trim, then $f(p)=2p_1$.
\end{proposition}

\proof
 The proof is by induction on $\max(p)$. 
 If $\max(p)=1$, then $p=(1)$ and the implication holds trivially ($|p\mkern1mu|>1$
 fails). Suppose now that $\max(p)>1$ and that the result holds for every trim element of $\ofs$ with smaller
 maximum entry. If $|p\mkern1mu|=1$, then again the implication holds trivially. Suppose that
 $n=|p\mkern1mu|>1$.  If $|R(p)|=1$, then $R(p)=(\gcdp{p})$, and thus
 $p=(\gcdp{p},2\gcdp{p})$, which is not trim. Thus $|R(p)|>1$. Suppose first that $R(p)$ is trim. Then by our 
 inductive hypothesis, $f(R(p))>\max(R(p))$. If $p_n\ge 2p_1$, then $\max(R(p))=p_n-p_1$ and thus 
 $f(p)=p_1+f(R(p))>p_1+p_n-p_1=p_n$. Now consider the case when $p_n<2p_1$. Then by Corollary \ref{frp}, 
 $f(p)=p_1+f(R(p))\ge p_1+p_1=2p_1>p_n$. Thus if $R(p)$ is
 trim, then $f(p)>p_n$. Suppose now that $R(p)$ is not trim. Then by Proposition \ref{ready for optimality} (ii),
 $f(R(p))=\max(R(p))$. If $p_n<2p_1$, then $\max(R(p))=p_1>p_n-p_1$, and then $f(p)=p_1+f(R(p))=2p_1>p_n$, as
 required. Suppose next that $p_n> 2p_1$. Then $\max(R(p))=p_n-p_1$, $|R(p)|=n$, and $R(p)\rest{n-1}=R(p\rest{n-1})$. 
 Since $R(p)$ is not trim, we have $\gcdp{R(p)}=\gcdp{R(p)\rest{|R(p)|-1}}$ and $\max(R(p))\ge f(R(p)\rest{|R(p)|-1})$, so
 $$
 \gcdp{p}=\gcdp{R(p)}=\gcdp{R(p)\rest{n-1}}=\gcdp{R(p\rest{n-1})}=\gcdp{p\rest{n-1}}
 $$
 and $p_n-p_1=\max(R(p))\ge f(R(p)\rest{n-1})=f(R(p\rest{n-1})$. Since $p$ is trim and $|p|>1$,
 this implies that $p_n<f(p\rest{n-1})=p_1+f(R(p\rest{n-1}))\le p_1+p_n-p_1=p_n$, which is impossible. Thus if $R(p)$ is 
 not trim, but $p$ is trim, it is not possible that $p_n>2p_1$.
 Finally, suppose that $p_n=2p_1$. Then $\gcdp{p}=\gcdp{p\rest{n-1}}$, so since $p$ is trim, we have $p_n<f(p\rest{n-1})$.
 Now $p_n-p_1=p_1$, so $R(p)=(p_2-p_1,\ldots,p_{n-1}-p_1,p_1)=R(p\rest{n-1})$, so $f(R(p))=f(R(p\rest{n-1})$, and
 thus 
 $$
 f(p)=p_1+f(R(p))=p_1+f(R(p\rest{n-1}))=f(p\rest{n-1})>p_n.
 $$
 As well, $\max(R(p))=p_1$, and so by Corollary \ref{not trim}, 
 $f(p)=p_1+f(R(p))=p_1+\max(R(p))=2p_1$. This completes the proof of the inductive step.
\endproof

 The following result gives an upper bound for $f$ that is reminiscent of the Fine-Wilf theorem.
 Later (see Proposition \ref{upper bound for fw(p)}), we shall
 establish a generalization of this which for $p$ trim with $|p\mkern1mu|\ge3$ offers a slightly improved upper bound for
 $fw$.
 
\begin{proposition}\label{upper bound for f}
 For $p\in\ofs$, $f(p)\le\min(p)+\max(p)-\gcdp{p}$.
\end{proposition}

\proof
 By Proposition \ref{reduced}, it suffices to prove that if $\gcdp{p}=1$, then $f(p)\le\min(p)+\max(p)-1$.
 The proof is by induction on $\max(p)$, with the base case $\max(p)=1$, so $p=(1)$ and $f(p)=1=\min(p)+\max(p)-1$.
 Suppose now that $\gcdp{p}=1$ and $\max(p)>1$ and the result holds for every element of $\ofs$ with smaller maximum entry. If
 $|p\mkern1mu|=1$, then $p=(\gcdp{p})=(1)$ and so $\max(p)=1$, which is not the case. Thus $n=|p\mkern1mu|>1$. Since $\gcdp{R(p)}=\gcdp{p}=1$,
 we may apply the induction hypothesis to $R(p)$ to obtain that $f(p)=p_1+f(R(p))\le p_1+ \min(R(p))+\max(R(p))-1$. We consider
 three cases. The first occurs when $p_1\le p_2-p_1$, in which case $\min(R(p))=p_1$,
 and $\max(R(p))=p_n-p_1$. We have $f(p)\le p_1+p_1+p_n-p_1-1=p_1+p_n-1$, as required. Next, suppose that
 $p_2-p_1<p_1\le p_n-p_1$. Then $\min(R(p))=p_2-p_1$, $\max(R(p))=p_n-p_1$, and we note that $p_2<2p_1$. We have
 $f(p)\le p_1+p_2-p_1+p_n-p_1-1=p_2-p_1+p_n-1<2p_1-p_1+p_n-1=p_1+p_n-1$. Finally, suppose that $p_n-p_1<p_1$, so that
 $\min(R(p))=p_2-p_1$ and $\max(R(p))=p_1$. Then $f(p)\le p_1+p_2-p_1+p_1-1=p_1+p_2-1\le p_1+p_n-1$. This completes the
 proof of the inductive step.
\endproof

\begin{corollary}\label{first upper bound fw}
 For $p\in\ofs$, $fw(p)\le \min(p)+\max(p)-\gcdp{p}$.
\end{corollary}

\proof
 By Proposition \ref{ready for optimality}, $fw(p)\le f(p)$, and by Proposition \ref{upper bound for f}, $f(p)\le \min(p)+\max(p)-\gcdp{p}$.
\endproof
 
\begin{theorem}[Fine-Wilf]\label{fine-wilf}
 Let $p\in\ofs$ be trim with $|p\mkern1mu|=2$. Then $fw(p)=p_1+p_2-\gcdp{p}$.
\end{theorem}

\proof
 Since $p$ is trim, $fw(p)=f(p)$, and by Proposition \ref{good lower bound for f}, $f(p)\ge p_1+p_2-\gcdp{p}$,
 while by Proposition \ref{upper bound for f}, $f(p)\le p_1+p_2-\gcdp{p}$.
\endproof

\section{The Fine-Wilfs graphs $G(p,k)$}

\begin{definition}
 Let $p\in \ofs$. For any $k\in \zplus$, $G(p,k)$ shall denote the simple graph with vertex set $\set 1,\ldots,k\endset$ and 
 edge set 
 $$
  \set \set i,j\endset\mid |i-j|=p_k\text{ for some $k=1,2,\ldots,|p\mkern1mu|$}\endset.
 $$
 The values $k=fw(p)$ and $k=fw(p)-1$ feature prominently in the development of the theory, and we shall
 let $G_p$ and $G_p'$ denote $G(p,fw(p))$ and $G(p,fw(p)-1)$, respectively.
\end{definition} 

Note that if $p,q\in\ofs$ and $q$ is the trimmed form of $p$, then $fw(p)=fw(q)$, and $G_p=G_q$, $G'_p=G'_q$.

Our first goal in this section is to establish that for $p\in\ofs$, the graph $G_p$ has exactly
$d=\gcdp{p}$ connected components, each isomorphic to $G_{p/d}$. For example,
for $p=(6,8,10)$, $fw(p)=12$, $\gcdp{p}=2$, and $G_p$ has two connected components, each isomorphic
to the connected graph $G_{(3,4,5)}$, where by Proposition \ref{reduced}, $fw(3,4,5)=fw(6,8,10)/2=6$.

{\centering
\table{c@{\hskip50pt}c}
$\vcenter{\xy /r20pt/:,
(0,0)*+{4}="4",
(1,0)*+{1}="1",
(2,0)*+{6}="6",
(3,0)*+{3}="3",
"1"+(0,-1)*+{5}="5",
"6"+(0,-1)*+{2}="2",
"4";"1"**\dir{-};"6"**\dir{-};"3"**\dir{-},
"1";"5"**\dir{-};"2"**\dir{-};"6"**\dir{-},
\endxy}$
&
$\vcenter{\xy /r20pt/:,
(0,0)*+{7}="7",
(1,0)*+{1}="1",
(2,0)*+{11}="11",
(3,0)*+{5}="5",
"1"+(0,-1)*+{9}="9",
"11"+(0,-1)*+{3}="3",
"7";"1"**\dir{-};"11"**\dir{-};"5"**\dir{-},
"1";"9"**\dir{-};"3"**\dir{-};"11"**\dir{-},
(4,0)+(0,0)*+{8}="8",
(4,0)+(1,0)*+{2}="2",
(4,0)+(2,0)*+{12}="12",
(4,0)+(3,0)*+{6}="6",
"2"+(0,-1)*+{10}="10",
"12"+(0,-1)*+{4}="4",
"8";"2"**\dir{-};"12"**\dir{-};"6"**\dir{-},
"2";"10"**\dir{-};"4"**\dir{-};"12"**\dir{-},
\endxy}$\\
\noalign{\vskip4pt}
$G_{(3,4,5)}$ & $G_{(6,8,10)}$
\endtable\par}

The second major objective of the section is then to establish that $G'_p$ has more than $d$ components.

\begin{definition}\label{component map}
 Let $p\in\ofs$ with $|p\mkern1mu|>1$. For any $k\ge 1$, the function $\alpha_{p,k}\from G(R(p),k)\to G(p,p_1+k)$ is
 defined by $\alpha_{p,k}(i)=p_1+i$.
\end{definition}

We note that in general, $\alpha_{p,k}$ is not a graph homomorphism. However, it does have the following
important property.

\begin{lemma}\label{alpha property}
 Let $p\in\ofs$ with $|p\mkern1mu|>1$, and let $k\ge 1$. If $i$ and $j$ belong to the same connected component of
 $G(R(p),k)$, then $\alpha_{p,k}(i)$ and $\alpha_{p,k}(j)$ belong to the same connected component of $G(p,p_1+k)$.
\end{lemma} 

\proof
 Let $i$ and $j$ be such that $1\le i<j\le k$ and $\set i,j\endset$ is an edge in $G(R(p),k)$.
 Suppose first of all that $j-i=p_1$. Then 
 $\set j,i\endset$ is an edge in $G(p,p_1+k)$, and since $\set i,\alpha_{p,k}(i)\endset$ and 
 $\set j,\alpha_{p,k}(j)\endset$ 
 are also edges in $G(p,p_1+k)$, it follows that $\alpha_{p,k}(i)$ and $\alpha_{p,k}(j)$ are connected in $G(p,p_1+k)$. 
 Now suppose that $j-i=p_t-p_1$ for some $t$ with $2\le t\le\max(p)$. Then $\set \alpha_{p,k}(j),i\endset$ is an edge in 
 $G(p,p_1+k)$. Since $\set i,\alpha_{p,k}(i)\endset$ is an edge in $G(p,p_1+k)$, it follows that $\alpha_{p,k}(i)$ and 
 $\alpha_{p,k}(j)$ belong to the same connected component of $G(p,p_1+k)$. 
\endproof

Thus if $C$ is a connected component of $G(R(p),k)$, then $\alpha_{p,k}(C)$
 is contained in a component of $G(p,p_1+k)$. 
 
We are now ready to demonstrate that for any $p\in\ofs$ with $\gcdp{p}=1$,  $\displaystyle G_{p/d}$
is connected. In fact, we have the following slightly stronger result.

\begin{proposition}\label{bigger than fw(p)}
 Let $p\in\ofs$ with $\gcdp{p}=1$, and let $k\ge fw(p)$. Then $G(p,k)$ is connected.
\end{proposition}

\proof
 The proof is by induction on $\max(p)$.
 The base case occurs when $\max(p)=1$ (since $\gcdp{(1)}=1$), and in this case, 
 $fw(p)=1$. Since for any $k\ge fw(p)=1$, $G(p,k)$ is a chain of length $k-1$, the result holds when $\max(p)=1$.
 Suppose now that $p\in\ofs$ has $\gcdp{p}=1$ and  $\max(p)>1$ and that the result holds for all elements of
 $\ofs$ of smaller maximum entry and greatest common divisor equal to 1. Note that $\gcdp{p}=1$ and $\max(p)>1$ imply
 that $n=|p\mkern1mu|>1$.
 
 \noindent Case 1: $p$ is not trim. Then $\gcdp{p}=\gcdp{p\rest{n-1}}$, $fw(p\rest{n-1})=fw(p)\le f(p)=p_n$, so 
 $G(p,fw(p))=G(p,fw(p\rest{n-1})=
 G(p\rest{n-1},fw(p\rest{n-1}))$. Since $\gcdp{p}=1$, we have $\gcdp{p\rest{n-1}}=1$ and so by hypothesis,  
 $G(p\rest{n-1},fw(p\rest{n-1}))$ is connected, which establishes that $G(p,fw(p))$ is connected.
 Suppose now that $k\ge fw(p)$. If $p_1>1$, then by Proposition \ref{fw lower bound}, $fw(p)\ge 2p_1>p_1$, while
 if $p_1=1$, then by Proposition \ref{min(p)=gcd(p)}, $fw(p)=p_1$. In any event, we have $fw(p)\ge p_1$.
 Thus for any vertex $i$ in $G(p,k)$ with $i>fw(p)\ge p_1$, we may write $i=qp_1+j$ where $1\le j\le p_1$, so that
 $j$ is a vertex of $G(p,fw(p))$ and $i$ is connected to $j$ in $G(p,k)$. Thus $G(p,k)$ is connected.
 
 \noindent Case 2: $p$ is trim. In this case, $fw(p)=f(p)$, so $G(p,fw(p))=
 G(p,f(p))=G(p,p_1+f(R(p)))$. Now, $1=\gcdp{p}=\gcdp{R(p)}$ and $\max(R(p))<\max(p)$, and since
 $f(R(p))\ge fw(R(p))$, it follows from the induction hypothesis that $G(R(p),f(R(p))$ is connected. Thus 
 $$
 \set p_1+1,p_1+2,\ldots,p_1+f(R(p))\endset,
 $$
 the image of $G(R(p),f(R(p)))$ under $\alpha_{R(p),f(R(p))}$, is a connected subset of $G(p,p_1+f(R(p)))=G(p,f(p))$.
 Now, for any $i$ with $1\le i\le p_1$, $i$ is connected to $i+p_1\le 2p_1\le fw(p)$ in $G(p,f(p))$,
 so $G(p,f(p))$ is connected. Since $G(p,f(p))=G(p,fw(p))$, $G(p,fw(p))$ is connected. Suppose now that $k\ge fw(p)$. 
 Then $G(p,fw(p))$ is a connected subgraph
 of $G(p,k)$. Since $p$ is trim and $|p\mkern1mu|>1$,  $fw(p)=f(p)\ge 2p_1$, so for $i$ in $G(p,k)$ with $i>fw(p)\ge 2p_1$, 
 we may write $i=qp_1+j$ where $1\le j\le p_1$. Thus
 $j$ is a vertex of $G(p,fw(p))$ and $i$ is connected to $j$ in $G(p,k)$, which proves that $G(p,k)$ is connected.
\endproof 

\begin{proposition}\label{connectedness in G(p,k)}
 Let $p\in \ofs$, $k\in\zplus$. If $i$ and $j$ belong to the same connected component of $G(p,k)$, then
 $\gcdp{p}$ divides $i-j$.
\end{proposition}

\proof
 It suffices to observe that if $i$ and $j$ are adjacent in $G(p,k)$, then $|i-j|=p_k$ for some $k$ with
 $1\le k\le |p\mkern1mu|$. Then since $\gcdp{p}$ divides $p_k$, it follows that $\gcdp{p}$ divides $i-j$.
\endproof

\begin{corollary}\label{lower bound} 
 For $p\in \ofs$ and $k\in \zplus$ with $k\ge \gcdp{p}$, $G(p,k)$ has at least $\gcdp{p}$ connected components.
\end{corollary}

\proof
 If $1\le i<j\le \gcdp{p}$, then it follows from Proposition \ref{connectedness in G(p,k)} that 
 $i$ and $j$ can't be in the same connected component of $G(p,k)$.
\endproof 

\begin{proposition}\label{components}
 Let $p\in\ofs$, and let $d=\gcdp{p}$. Then for each $i=1,2,\ldots,d$, the map $\gamma_i\from G_{p/d}\to
 G_p$ defined by $\gamma_i(j)=i+(j-1)d$ for $1\le j\le fw(p/d)$ is an injective graph homomorphism which 
 is an isomorphism from $G_{p/d}$ onto the subgraph $\gamma_i(G_{p/d})$ of $G_p$.
 Moreover, $G_p$ has exactly $d$ components, the images of $\gamma_i$, $i=1,2,\ldots,d$; that is, the congruence
 classes of the interval $\{\,1,2,\ldots,fw(p)\,\}$.
\end{proposition}

\proof
 It is immediate from Proposition \ref{connectedness in G(p,k)} that each component of $G_p$ is contained
 in the image of $\gamma_i$ for some $i$ with $1\le i\le d$, and by Proposition \ref{bigger than fw(p)}
 $G_{p/d}$ is connected. It remains only to prove that for each
 such $i$, $\gamma_i$ is a graph homomorphism, injectivity being obvious. Let $j$, $k$ be vertices
 of $G_{p/d}$. Since $|\gamma_i(j)-\gamma_i(k)|=|(j-1)d-(k-1)d|=
 |j-k|d$, it follows that $|i-j|=p_t/d$ if and only if $|\gamma_i(j)-\gamma_i(k)|=p_t$.
 Thus $\gamma_i$ is a graph isomorphism from $G_{p/d}$ onto the subgraph $\gamma_i(G_{p/d})$ of $G_p$.
\endproof

\begin{corollary}\label{component count for G(p,k)}
 Let $p\in\ofs$ and let $k\ge fw(p)$. Then $G(p,k)$ has exactly $d=\gcdp{p}$ components, the congruence classes
 of the interval $\{\,1,2,\ldots,k\,\}$ modulo $d$.
\end{corollary}

\proof
 Since $G_p$ is a subgraph of $G(p,k)$, and by Proposition \ref{components} , $G_p$ has exactly
 $\gcdp{p}$ components, it suffices to prove that each vertex $i$ of $G(p,k)$ is connected
 to a vertex in the subgraph $G_p$. But $i=qp_1+j$, where $1\le j\le p_1$, and $j$ is
 a vertex of $G_p$, with $i$ connected to $j$ by a path of length $q$ in $G(p,k)$.
\endproof 

\begin{proposition}\label{alpha injective on components}
 Let $p\in\ofs$ have $|p\mkern1mu|>1$, and let $k\ge 1$. 
 Then $\alpha_{p,k}$ induces an injective map from the set of components of $G(R(p),k)$ into the set of 
 components of $G(p,p_1+k)$.
\end{proposition}

\proof
 Suppose not, and of all pairs of disconnected elements of $G(R(p),k)$ whose images are
 connected in $G(p,p_1+k)$, choose two, say $i<j$, whose shortest path joining their images in $G(p,p_1+k)$ has 
 least length. Let $i_0=i+p_1,i_1,\ldots,i_n=j+p_1$ denote a shortest path in $G(p,p_1+k)$ from $i+p_1$ to $j+p_1$.
 Let $m\ge0$ be the maximum index such that for all $t$ with $0\le t\le m$, $i_t> p_1$. Suppose that $i_1>p_1$.
 If $i_0>i_1$, then $i_0-i_1=p_t$ for some $t$, and thus $i-i_1=p_t-p_1$. But then $i\ge i_1>p_1$, and thus
 $(i-p_1)-(i_1-p_1)=p_t-p_1$, which means that $\set i-p_1,i_1-p_1\endset$ is an edge in $G(R(p),k)$. Since
 $\set i,i-p_1\endset$ would also be an edge in $G(R(p),k)$, this would imply that $i$ and $i_1-p_1$ are in the
 same component of $G(R(p),k)$, and thus $i_1-p_1$ and $j$ are in different components of $G(R(p),k)$. Since
 $(i_1-p_1)+p_1=i_1$ and $j+p_1$ are connected in $G(p,k+p_1)$ by a path of length $n-1$, we have a contradiction
 to the minimality of $n$. Thus if $i_1>p_1$, it must be that $i_0<i_1$ (since $i_0=i_1$ is not possible). In this
 case, $i_1-i_0=p_t$ for some $t$, and so $(i_1-p_1)-i=p_t$. Thus $i_1-p_1> p_t\ge p_1$, and so $i_1-2p_1>0$,
 from which we obtain that $(i_1-2p_1)-i=p_t-p_1\ge0$. If $(i_1-2p_1)-i=0$, then $i_1-2p_1=i$, otherwise 
 $\set i_1-2p_1,i\endset$ is an edge in $G(R(p),k)$.
 Since $i_1\le k+p_1$, we have $i_1-p_1\le k$ and thus $\set i_1-p_1,i_1-2p_1\endset$ is an edge in $G(R(p),k)$ as
 well. Thus $i_1-p_1$ lies in the same component of $G(R(p),k)$ as does $i$, which means that $i_1-p_1$ and $j$
 lie in different components of $G(R(p),k)$, again contradicting the minimality of $n$. Thus $i_1\le p_1$.
 We consider two cases: $i_1\le k$, and $i_1>k$. Suppose first that $i_1\le k$. Since $i_1\le p_1$, we have
 $i+p_1>i_1$, and so $(i+p_1)-i_1=p_t$ for some $t$. But then $i-i_1=p_t-p_1\ge0$, and since $i,i_1\le k$, with $i+p_1=i_0\ne i_1$,
 it follows
 that either $i=i_1$ or else $\set i,i_1\endset$ is an edge in $G(R(p),k)$. Now since $i_n=j+p_1>p_1$, while $i_1\le p_1$,
 we conclude that $n\ge 2$. If $i_2\le p_1$, then
 $|i_1-i_2|<p_1$, contradicting the fact that $|i_1-i_2|=p_r$ for some $r$. Thus $i_2> p_1\ge i_1$, and so $i_2-i_1=p_r$
 for some $r$. Now $p_1<i_2\le k+p_1$, so $0<i_2-p_1\le k$ and $(i_2-p_1)-i_1=p_r-p_1$, so $\set (i_2-p_1),i_1\endset$
 is an edge in $G(R(p),k)$. Thus $i_2-p_1$ and $i$ lie in the same component of $G(R(p),k)$, contradicting the
 minimality of $n$. Thus $i_1\le k$ is impossible, which means that we must have $i_1>k$. But then $i\le k<i_1$.
 However, since $i_1<i+p_1$ and $\set i_1,i+p_1\endset$ is an edge in $G(R(p),k)$, we have $(i+p_1)-i_1=p_r$
 for some $r$. But then $i-i_1=p_r-p_1\ge 0$ and thus $i\ge i_1$, impossible. It follows therefore that the
 map on components that is induced by $\alpha_{p,k}$ is injective.
\endproof

\begin{proposition}\label{optimality}
 Let $p\in\ofs$. Then $\min(p)\ne \gcdp{p}$ implies that $G'_p$ has more than $\gcdp{p}$ components.
\end{proposition}

\proof
 As usual, the proof is by induction on $\max(p)$, and the result is trivially true for $\max(p)=1$.
 Suppose then that $p\in\ofs $ has $\max(p)>1$ and the result holds for all elements of $\ofs$ with
 smaller maximum entry. Let $d=\gcdp{p}$ and $n=|p\mkern1mu|$, and suppose that $p_1\ne d$. 
 We first consider the case when $p$ is not trim. In this case, it follows from Corollary \ref{not trim} that
 $p_n=f(p)\ge fw(p\rest{n-1})=
 fw(p)$, so $G_p'=G_{p\rest{n-1}}'$. The inductive hypothesis can therefore be applied to
 conclude that $G_{p\rest{n-1}}'$ has more than $\gcdp{d}=\gcdp{p\rest{n-1}}$ components.
 
 Assume now that $p$ is trim.
 
 Suppose first that $R(p)$ is trim. If $\min(R(p)) = \gcdp{R(p)}$, then by Lemma \ref{benefits of trim}, $|R(p)\mkern.8mu|=1$
 and so $R(p)=(d)$ since $d=\gcdp{p}=\gcdp{R(p)}$. But then $p=(d,2d)$, so $\gcdp{p}=\min(p)$, which is not
 the case. Thus $\min(R(p)) \ne \gcdp{R(p)}$, and so by the inductive
 hypothesis, $G_{R(p)}'$ has more than $\gcdp{p}=\gcdp{R(p)}$ components.
 But $fw(R(p))=f(R(p))=f(p)-p_1=fw(p)-p_1$, and so $G_{R(p)}'=G(R(p),fw(p)-p_1-1)$ has more 
 than $\gcdp{R(p)}=\gcdp{p}$ components. By Proposition \ref{alpha injective on components}, 
 $\alpha_{p,fw(p)-1-p_1}\from G(R(p),fw(p)-1-p_1)\to G(p,fw(p)-1)$ induces an injective function
 on components, and so it follows that $G_p'=G(p,fw(p)-1)$ has more than $\gcdp{p}$ components.
 
 Suppose now that $R(p)$ is not trim. Then by Proposition \ref{lower bound for fw},
 $f(p)=2p_1$. Since $p$ is trim, we have $fw(p)=f(p)=2p_1$, and so it follows that 
 $G_p'=G(p,fw(p)-1)=G(p,2p_1-1)$ has $\set p_1\endset$ as a
 component. By Proposition \ref{alpha injective on components}, the map 
 $$
 \alpha_{p,p_1-1}\from
 G(R(p),p_1-1)=G(R(p),fw(p)-p_1-1)\to G(p,fw(p)-1)=G_pf
 $$
 induces an
 injective map on components. Since $p_1$ is not in the image of $\alpha_{p,p_1-1}$,  $G_p'$
 has at least one more component than $G(R(p),p_1-1)$.
 We have $p_1=f(R(p))\ge fw(R(p))$. If $fw(R(p))<p_1$, then $G(R(p),p_1-1)$ has $\gcdp{R(p)}=\gcdp{p}$
 components, which then implies that $G_p'$ has more than $\gcdp{p}$ components. 
 It remains to consider the case when $fw(R(p))=p_1$. Then $G(R(p),p_1-1)=G(R(p),fw(R(p))-1)=G_{R(p)}'$. 
 Suppose that $\min(R(p)) = \gcdp{R(p)}$. Then by Proposition \ref{min(p)=gcd(p)}, 
 $p_1=fw(R(p))=\min(R(p))=\gcdp{R(p)}=\gcdp{p}$, which is not the case.
 Thus $\min(R(p)) \ne \gcdp{R(p)}$ and we may therefore apply the induction hyposthesis to $G_{R(p)}'=
 G(R(p),p_1-1)$, to conclude that $G(R(p),p_1-1)$ has more than $\gcdp{R(p)}=\gcdp{p}$ components, which then
 implies that $G(p,fw(p)-1)$ has more than $\gcdp{p}$ components. This completes the proof of the inductive step.
\endproof

Note: Let $m(p)=fw(p)-1$. The relationship between $m(p)$ and $m(p/d)$, $d=\gcdp{p}$ is quite straightforward.
By Proposition \ref{reduced}, $fw(p)=dfw(p/d)$, so $m(p)=fw(p)-1=dfw(p/d)-1=d(fw(p/d)-1+1)-1=d(m(p/d)+1)-1$.

\begin{proposition}\label{symmetry}
 Let $p\in\ofs$. For each positive integer $k$, the map $\tau_{p,k}\from G(p,k)\to
 G(p,k)$ given by $\tau_{p,k}(i)=  k-(i-1)$ is a graph automorphism of order 2.
\end{proposition}

\proof
 For $i\in\z$, $1\le i\le k$ if and only if $-1\ge -i\ge -k$, which in turn
 holds if and only if $k\ge k-i+1\ge 1$, so $i$ is a vertex of $G(p,k)$ if and only if
 $k-(i-1)$ is a vertex of $G(p,k)$. Thus $\tau_{p,k}$ is a function from the vertex set of
 $G(p,k)$ to itself, evidently of order 2. Furthermore, for $1\le i<j\le k$, $j-i=p_t$ for some
 $t$ if and only if $\tau_{p,k}(i)-\tau_{p,k}(j)=(k-i+1)-(k-j+1)=j-i=p_t$, and so $\set i,j\endset$ is a edge in $G(p,k)$ if
 and only if $\set \tau_{p,k}(i),\tau_{p,k}(j)\endset$ is an edge in $G(p,k)$.
\endproof

Of course, since $\tau_{p,k}$ is a graph automorphism of $G(p,k)$, it induces a permutation of the set
of components of $G(p,k)$, and in general, this permutation is nontrivial. For example, $p=(8,12,14)$ is a trim
sequence, so $fw(p)=f(p)=16$, and $\tau_{p,fw(p)}\from G_p\to G_p$ induces a nontrivial permutation
on the set consisting of the two components of $G_p$. This may be quickly verified by observing that since
$\tau_{p,fw(p)}(i)=17-i$ for each $i$, we have in particular that $\tau_{p,fw(p)}(6)=17-6=11$. 

{\centering
\table{c@{\hskip50pt}c}
$\vcenter{\xy /r20pt/:,
(0,0)*+{6}="6",
(1,0)*+{14}="14",
(2,0)*+{2}="2",
(3,0)*+{16}="16",
(4,0)*+{4}="4",
(5,0)*+{12}="12",
"2"+(0,-1)*+{10}="10",
"16"+(0,-1)*+{8}="8",
"6";"14"**\dir{-};"2"**\dir{-};"16"**\dir{-};"4"**\dir{-};"12"**\dir{-},
"2";"10"**\dir{-},
"16";"8"**\dir{-},
\endxy}$
&
$\vcenter{\xy /r20pt/:,
(0,0)*+{11}="6",
(1,0)*+{3}="14",
(2,0)*+{15}="2",
(3,0)*+{1}="16",
(4,0)*+{13}="4",
(5,0)*+{5}="12",
"2"+(0,-1)*+{7}="10",
"16"+(0,-1)*+{9}="8",
"6";"14"**\dir{-};"2"**\dir{-};"16"**\dir{-};"4"**\dir{-};"12"**\dir{-},
"2";"10"**\dir{-},
"16";"8"**\dir{-},
\endxy}$
\\
\noalign{\vskip4pt}
\multicolumn{2}{c}{$G((8,12,14),16)$}
\endtable\par}

More generally, we may easily describe the fixed points of the permutation that $\tau_{p,fw(p)}$ induces
on the set of components of
$G_p$. By Proposition \ref{components}, $i$ and $j$ are in the same component of $G_p$ if 
and only if $i\cong j\mod{d}$, where $d=\gcdp{p}$. Since $d$ divides $fw(p)$, it follows that for any $i$, we have
$i-\tau_{p,fw(p)}(i)=2i-1-fw(p)$, so $\tau_{p,fw(p)}$ fixes the component of $i$ (setwise) if
and only if $d$ divides $2i-1$. Thus in general, the permutation on the set of components of $G_p$
that is induced by $\tau_{p,fw(p)}$ will have few fixed points. The following result thus comes
as a bit of a surprise.

\begin{proposition}\label{palindrome}
 Let $p\in\ofs$. If $\gcdp{p}<p_1$, then $\tau_{p,fw(p)-1}(C)=C$ for each component $C$ of $G'_p$.
\end{proposition}

\proof
 We are to prove that for each $i$ with $1\le i\le fw(p)-1$, $i$ and $\tau_{p,fw(p)-1}(i)=fw(p)-1-(i-1)=fw(p)-i$
 are connected
 in $G'_p$, and as usual, we shall use induction on $\max(p)$. If $\max(p)=1$, then $p=(1)$ and $\gcdp{p}=1$, so 
 the result holds vacuously. 
 Suppose now that
 $\max(p)>1$, and that the result holds for all elements of $\ofs$ with maximum entry less than
 $\max(p)$. If $|p\mkern1mu|=1$, then again, we have $\gcdp{p}=p_1$ and so the result holds vacuously.
 Suppose that $|p\mkern1mu|>1$. Suppose first that $p$ is not trim, and let $q$ denote the trimmed form of $p$.
 Since $fw(p)=fw(q)$, $G_p'=G_q'$, $\tau_{p,fw(p)-1}=]tau_{q,fw(q)-1}$, and $\max(q)<\max(p)$, we may apply the
 induction hypothesis to $q$, so the result holds for $q$ and thus for $p$. We may therefore assume that $p$ is
 trim, so $fw(p)=f(p)$. Let $i$ be such that $1\le i\le f(p)-1$.
 Note that if $1\le i\le p_1$, then $f(p)-1\ge f(p)-i\ge f(p)-p_1$,  and by Lemma \ref{f bounds}, 
 $f(p)\ge 2p_1$, so $f(p)-1\ge f(p)-i\ge p_1$. Since $\tau_{p,f(p)-1}$ is an automorphism of order 2, we may 
 assume that $i\ge p_1$.
 
 We apply the induction hypothesis to $R(p)$ to conclude that $i-p_1$ and $fw(R(p))-(i-p_1)$ are connected in $G'_{R(p)}$.
 We claim that $i-p_1$ and $f(R(p))-(i-p_1)$ are connected in $G(R(p),f(R(p))-1)$.
 If $fw(R(p))=f(R(p)$, there is nothing to show, while if $fw(R(p))<f(R(p))$, then by Proposition \ref{component count for G(p,k)},
 the components of $G(R(p)),f(R(p))-1)$ are precisely the congruence classes of $\{\,1,2,\ldots,f(R(p))-1\,\}$ modulo $d=\gcd{R(p)}=
 \gcdp{p}$, and since $f(R(p))\cong 0\cong f(p)\mod{d}$, it follows that $i-p_1\cong fw(R(p))-(i-p_1)\cong f(R(p))-(i-p_1)\mod{d}$ and
 so $i-p_1$ and $f(R(p))-(i-p_1)$ are connected in $G(R(p),f(R(p))-1)$, as required.
 But then by Proposition \ref{alpha property}, $i$ and $p_1+f(R(p))-i+p_1=f(p)-(i-p_1)=fw(p)-(i-p_1)$ are connected in $G(p,p_1+f(R(p))-1)=G_p'$. 
 Since $f(p)-i$ and $f(p)-i+p_1$ are connected in $G'_p$ as well, we finally obtain that $i$ and $fw(p)-i$ are connected in $G'_p$,
 as required. This completes the proof of the inductive step, and so the result
 follows.
\endproof

We have seen that for $k>fw(p)-1$, in particular for $k=fw(p)$, that the permutation on the set of components of
$G(p,k)$ by $\tau_{p,k}$ may have many fixed points, and the same holds true for $k<fw(p)-1$.
For example, if $p=(3,5)$, then $fw(p)=7$, so $fw(p)-1=6$. Let us consider $k=5$. The components
of $G(p,5)$ are $\set 1,4\endset$, $\set 2,5\endset$, and $\set 3\endset$, and we have
$\tau_p(\set 1,4\endset=\set 5-(1-1),5-(4-1)\endset=\set 2,5\endset$, so 
$\tau_p(\set 2,5\endset)=\set 1,4\endset$, while $\tau_p(\set 3\endset)=\set 3\endset)$.

{\centering
\table{c@{\hskip50pt}c}
$\vcenter{\xy /r20pt/:,
(0,0)*+{1}="1",
(1,0)*+{4}="4",
(2,0)*+{3}="3",
(3,0)*+{2}="2",
(4,0)*+{5}="5",
"1";"4"**\dir{-},
"2";"5"**\dir{-},
\endxy}$
&
$\vcenter{\xy /r20pt/:,
(0,0)*+{3}="3",
(1,0)*+{6}="6",
(2,0)*+{1}="1",
(3,0)*+{4}="4",
(4,0)*+{2}="2",
(5,0)*+{5}="5",
"3";"6"**\dir{-},
"6";"1"**\dir{-},
"1";"4"**\dir{-},
"2";"5"**\dir{-},
\endxy}$
\\
\noalign{\vskip4pt}
$k=5$, $\tau_p(i)=6-i$ & $k=6$, $\tau_p(i)=7-i$
\endtable\par}

We remark that in general, $|\aut(G_p)|>2$. For example, if $p=(8,12,18,19)$, then $\aut(G_p)\simeq (C_2^2\times C_2^2)\rtimes C_2$,
where the action of $C_2$ on $C_2\times C_2)$ is to swap coordinates, while $\aut(G_p')\simeq D_5\times C_2\times ( (C_2\times C_2)\rtimes C_2)$,
where again, the action of $C_2$ is to swap coordinates.

\begin{example} Let $p=(8,12,18,19)$. Then

\vskip6pt
\hbox to \hsize{\hfil
\table{@{\hskip-6.2pt}c@{\hskip6pt}c@{}}
\def\nodenames{\ifcase\xypolynode\or3\or16\or1\or14\or2\or15\fi}
$\vcenter{\xy /r19pt/:,
(0,0)*+{\hbox{\scriptsize$1$}}="1",
(-1.5,0)*+{\hbox{\scriptsize$19$}}="19",
(-2.5,0)="o"+a(72)*+{\hbox{\scriptsize$11$}}="11",
"o"+a(144)*+{\hbox{\scriptsize$3$}}="3",
"o"+a(216)*+{\hbox{\scriptsize$15$}}="15",
"o"+a(288)*+{\hbox{\scriptsize$7$}}="7",
"1";"19"**\dir{-},
"19";"11"**\dir{-},
"11";"3"**\dir{-},
"3";"15"**\dir{-},
"15";"7"**\dir{-},
"7";"19"**\dir{-},
(0,-1)="v"+a(162)*+{\hbox{\scriptsize$9$}}="9",
"v"+a(234)*+{\hbox{\scriptsize$17$}}="17",
"v"+a(306)*+{\hbox{\scriptsize$5$}}="5",
"v"+a(18)*+{\hbox{\scriptsize$13$}}="13",
"1";"9"**\dir{-},
"9";"17"**\dir{-},
"17";"5"**\dir{-},
"5";"13"**\dir{-},
"13";"1"**\dir{-},
(2,0)*+{\hbox{\scriptsize$20$}}="20",
"1";"20"**\dir{-},
(2,1)="w"+a(54)*+{\hbox{\scriptsize$4$}}="4",
"w"+a(135)*+{\hbox{\scriptsize$16$}}="16",
"w"+a(198)*+{\hbox{\scriptsize$8$}}="8",
"w"+a(-18)*+{\hbox{\scriptsize$12$}}="12",
"20";"12"**\dir{-},
"12";"4"**\dir{-},
"4";"16"**\dir{-},
"16";"8"**\dir{-},
"8";"20"**\dir{-},
(3.5,0)*+{\hbox{\scriptsize$2$}}="2",
"20";"2"**\dir{-},
(4.5,0)="o"+a(108)*+{\hbox{\scriptsize$14$}}="14",
"o"+a(36)*+{\hbox{\scriptsize$6$}}="6",
"o"+a(-36)*+{\hbox{\scriptsize$18$}}="18",
"o"+a(-108)*+{\hbox{\scriptsize$10$}}="10",
"2";"14"**\dir{-},
"14";"6"**\dir{-},
"6";"18"**\dir{-},
"18";"10"**\dir{-},
"10";"2"**\dir{-},
\endxy}$
&
$\vcenter{\xy /r19pt/:,
(0,0)*+{\hbox{\scriptsize$1$}}="1",
(-1.5,0)*+{\hbox{\scriptsize$19$}}="19",
(-2.5,0)="o"+a(72)*+{\hbox{\scriptsize$11$}}="11",
"o"+a(144)*+{\hbox{\scriptsize$3$}}="3",
"o"+a(216)*+{\hbox{\scriptsize$15$}}="15",
"o"+a(288)*+{\hbox{\scriptsize$7$}}="7",
"1";"19"**\dir{-},
"19";"11"**\dir{-},
"11";"3"**\dir{-},
"3";"15"**\dir{-},
"15";"7"**\dir{-},
"7";"19"**\dir{-},
(0,-1)="v"+a(162)*+{\hbox{\scriptsize$9$}}="9",
"v"+a(234)*+{\hbox{\scriptsize$17$}}="17",
"v"+a(306)*+{\hbox{\scriptsize$5$}}="5",
"v"+a(18)*+{\hbox{\scriptsize$13$}}="13",
"1";"9"**\dir{-},
"9";"17"**\dir{-},
"17";"5"**\dir{-},
"5";"13"**\dir{-},
"13";"1"**\dir{-},
(1.2,1)="w"+a(54)*+{\hbox{\scriptsize$4$}}="4",
"w"+a(135)*+{\hbox{\scriptsize$16$}}="16",
"w"+a(198)*+{\hbox{\scriptsize$8$}}="8",
"w"+a(-18)*+{\hbox{\scriptsize$12$}}="12",
"12";"4"**\dir{-},
"4";"16"**\dir{-},
"16";"8"**\dir{-},
(2.5,0)*+{\hbox{\scriptsize$2$}}="2",
(3.5,0)="o"+a(108)*+{\hbox{\scriptsize$14$}}="14",
"o"+a(36)*+{\hbox{\scriptsize$6$}}="6",
"o"+a(-36)*+{\hbox{\scriptsize$18$}}="18",
"o"+a(-108)*+{\hbox{\scriptsize$10$}}="10",
"2";"14"**\dir{-},
"14";"6"**\dir{-},
"6";"18"**\dir{-},
"18";"10"**\dir{-},
"10";"2"**\dir{-},
\endxy}$\\
\noalign{\vskip6pt}
 $\aut(G_p)\simeq (C_2^2\times C_2^2)\rtimes C_2$ &  $\aut(G'_p)\simeq
  D_5\mkern-5mu\times\mkern-5mu C_2\mkern-5mu\times\mkern-5mu ((C_2\mkern-5mu\times\mkern-5mu C_2)\mkern-5mu\rtimes\mkern-3mu C_2)$.
\endtable\hfil}

\end{example}

\begin{proposition}\label{first reduction}
 Let $p\in\ofs$ be such that there exist $i$ and $j$ with $1\le i<j\le |p\mkern1mu|$ and $p_j$ is a multiple 
 of $p_i$, and let $q\in\ofs$ denote the sequence obtained by deleting $p_j$ from $p$.
 Then 
 \begin{enumerate}
  \item For every positive integer $k$, the components of $G(p,k)$ are identical to those
  of $G(q,k)$. In particular, $\component{G(p,k)}=\component{G(q,k)}$.
  \item $fw(p)=fw(q)$.
  \item $p_j\ge fw(p_{j-1})$, and $\gcdp{p\rest{j}}=\gcdp{p\rest{j-1}}$.
 \end{enumerate} 
\end{proposition}

\proof
 Suppose that $p_j=tp_i$.
 Let $1\le r<s\le k$, and suppose that $r$ and $s$ are joined by a walk in $G(p,k)$. If no
 edge in the walk is determined by $p_j$, then the walk is a walk from $r$ to $s$ in $G(q,k)$.
 Otherwise, there is at least one edge in the walk that is determined by $p_j$. If
 $\set a,b\endset$ is an edge in the walk, where $a>b$ and $a-b=p_j$, then we may replace
 $\set a,b\endset$ by the path $b, b+p_i, b+2p_i,\ldots,b+tp_i=b+p_j=a$. Apply this procedure to all
 edges in the walk that are determined by $p_j$ to obtain a walk from $r$ to $s$ that does
 not use any edge determined by $p_j$; that is, a walk from $r$ to $s$ in $G(q,k)$.
 Thus any component of $G(p,k)$ is contained in a component of $G(q,k)$. Since $G(q,k)$
 is a spanning subgraph of $G(p,k)$, each component of $G(q,k)$ is contained in a component
 of $G(p,k)$. Thus the components of $G(p,k)$ are identical to those of $G(q,k)$.
 
 Now since $\gcdp{p}=\gcdp{q}$, it follows from Corollary \ref{component count for G(p,k)}  and 
 Proposition \ref{optimality} that if $\gcdp{p}\ne p_1$, then 
 \begin{align*}
    fw(p)&=\min\set k\mid \component{G(p,k)}=\gcdp{p}\endset\\
    &=\min\set k\mid \component{G(q,k)}=\gcdp{p}\endset=fw(q),
 \end{align*}
 while if $\gcdp{p}=p_1$, then by Proposition \ref{min(p)=gcd(p)}, we have $fw(p)=p_1=fw(q)$.
 
 For (iii), we may assume that $\gcdp{p}=1$ and that $j=|p\mkern1mu|$. By Proposition \ref{optimality}, either $p_1=\min(p)=\gcdp{p}$ or else
 $G'_p$ is not connected. If $p_1=\gcdp{p}$, then since by (ii), $fw(p\rest{j-1})=fw(p)$, we have $fw(p\rest{j-1})=p_1<p_j$.
 Otherwise, $G'_p=G_p-fw(p)$ is not connected, while $G_p$ is connected, and so $fw(p)$ is a cut-vertex of $G_p$.
 Suppose that $p_j<fw(p\rest{j-1})$, so that by (ii), $p_j<fw(p)$, and thus $1+p_j\le fw(p)$. We have
 $p_j=cp_i$ for some integer $c$ and some $p_i$ in $p$, so there is a path from $1$ to $1+p_j$ in 
 $G_p$ consisting of $c$ edges determined by $p_1$. But then we may follow the edge from $1+p_j$ back to 1
 that is determined by $p_j$, and so $p_j$ belongs to a cycle in $G_p$. Since $p_j\le fw(p)$, there exists $t\ge0$
 such that $t+p_j=fw(p)$, and then $1+t, 1+t+p_i,\ldots,1+t+cp_i=1+t+p_n=fw(p),1+t$ is a cycle in $G_p$
 through $fw(p)$, which contradicts the fact that $fw(p)$ is a cut-vertex of $G_p$. Thus $p_j\ge fw(p\rest{j-1})$.
 Finally, since $p_j$ is a multiple of $p_i$, it is immediate that $\gcdp{p\rest{j}}=\gcdp{p\rest{j-1}}$.
\endproof 

The preceding result suggests that it will be convenient to introduce notation for those sequences with
no entry a multiple of another. 

\begin{definition}\label{reduced def}
 For $p\in\ofs$, the sequence obtained from $p$ by iterated deletion of any term that is a multiple of another shall
 be called the reduced form of $p$. Furthermore, we shall say that
 $p$ is reduced if $p$ is equal to its reduced form. 
\end{definition}

\begin{corollary}\label{reduced form has equal fw}
 If $p\in\ofs$ and $q$ is the reduced form of $p$, then $fw(p)=fw(q)$.
\end{corollary}

\proof
This is immediate from Proposition \ref{first reduction}.
\endproof

\begin{lemma}\label{implies connected}
 Let $p\in\ofs$, and let $k$ be a positive integer. If the interval $\set 1,2,\ldots,p_1\endset$ is contained within
 a component of $G(p,k)$, then $G(p,k)$ is connected.
\end{lemma}

\proof
 If $p_1=1$, then $G(p,k)$ is a complete graph, hence connected. Suppose that $p_1>1$, in which case the fact that
 1 and 2 are connected in $G(p,k)$ implies that $k>p_1$.
 For each $i$ with $p_1<i\le k$, there exist positive integers
 $m$ and $j$ such that $i=mp_1+j$ and $1\le j\le p_1$. Thus there is a path of length $m$ from
 $i$ to $j$ in $G(p,k)$. It follows now that $G(p,k)$ is connected.
\endproof

\begin{definition}\label{def of redundant}
 Let $p\in\ofs$. For $j$ such that $2\le j\le |p\mkern1mu|$, $p_j$ is said to be {\it redundant} in $p$ if
 $\gcdp{p\rest{j}}=\gcdp{p\rest{j-1}}$ and $p_j\ge f(p\rest{j-1})$.
\end{definition}
 
\begin{proposition}\label{second reduction}
 Let $p\in\ofs$ and $j$ be such that $p_j$ is redundant in $p$.
 If $q\in\ofs$ denotes the sequence obtained by deleting $p_j$ from $p$,
 then $fw(p)=fw(q)$.
\end{proposition}

\proof
 First, note that $\gcdp{p}=\gcdp{q}$. Let $n=|p\mkern1mu|$. If $p_1=\gcdp{p}$, then $fw(p)=p_1=fw(q)$, so we may
 suppose that $p_1\ne \gcdp{p}$. We may further assume that $\gcdp{p}=1$, and that $p$ is trim. Suppose 
 that $fw(p)< fw(q)$. By Proposition \ref{optimality}, $G(q,fw(p))$ is not connected, and by 
 Proposition \ref{lower bound for fw}, $p_j<fw(p)$. Thus $fw(p\rest{j-1})\le p_j<fw(p)$, and so $G(p\rest{j-1},fw(p\rest{j-1}))$
 is a connected subgraph of $G(q,fw(p))$. By Proposition \ref{fw lower bound} and the fact that $j-1\ge 2$ (since
 $j>1$ and if $j=2$, then from $\gcdp{p\rest{j}}=\gcdp{p\rest{j-1}}$ we would have $\gcdp{p}=p_1$), 
 we have $2p_1\le fw(p\rest{j-1})$, and so the interval $\set 1,2,\ldots,p_1\endset$ is contained within
 a component of $G(q,fw(p))$. By Lemma \ref{implies connected}, this implies that $G(q,fw(p))$ is
 connected, which is not the case. Therefore, it must
 be that $fw(p)\ge fw(q)$. Suppose that $fw(p)>fw(q)$. Then $G(q,fw(q))$ is a connected subgraph of
 $G'_p$. By Proposition \ref{fw lower bound} applied to $q$, $fw(q)\ge 2p_1$.
 Thus the interval $\set 1,2,\ldots,2p_1\endset$ is contained with a component of $G'_p$ and so by
 Lemma \ref{implies connected}, $G'_p$ is connected. But this is not the case, so we conclude that
 $fw(p)=fw(q)$, as required.
\endproof 

 As a consequence of Proposition \ref{second reduction}, we see that given $p\in \ofs$, we may iteratively
 remove redundant entries without regard to the order of removal to end up with a sequence with no
 redundant entries. More precisely, if we construct a list of elements of $\ofs$ with first entry
 $p$, and each subsequent entry obtained by selecting and removing a redundant element from the current entry in the list,
 then the last entry in the list will equal the element of $\ofs$ that is obtained from $p$ by identifying all redundant elements in $p$
 and removing them all at the same time. 
 
\begin{definition}\label{def of rp}
 For $p\in\ofs$, let $r(p)$ equal the number of redundant entries in $p$, and let $\hat{p}$ denote the 
 element of $\ofs$ that is obtained by deleting all $r(p)$ redundant entries in $p$, so $r(p)=|p\mkern1mu|-|\hat{p}|$.
 $\hat{p}$ shall be called the totally reduced form of $p$, and we shall say that $p$ is totally reduced
 if $p=\hat{p}$.
\end{definition}
 
\begin{corollary}\label{fw of totally reduced form}
 Let $p\in\ofs$. Then $fw(\hat{p})=fw(p)$.
\end{corollary}

\proof
 This follows from Proposition \ref{second reduction} and the fact that $\hat{p}$ can be formed from $p$ by
 $r(p)$ iterations of the process of selecting and removing a redundant entry.
\endproof

Note that if $p\in\ofs$ is totally reduced, then for every $j$ with $1\le j\le |p\mkern1mu|$, $p\rest{j}$ is both reduced 
and trim. 
 
We are now in a position to give an upper bound for $fw(p)$ that is an improvement over that given in 
Proposition \ref{upper bound for f} (provided that $r(p)<|p|-1$, its maximum possible value).

\begin{proposition}\label{upper bound for fw(p)}
 For each $p\in\ofs$,
 $$
  fw(p)\le \min(p)+\max(p)-\gcdp{p}(|p\mkern1mu|-1-r(p)).
 $$ 
\end{proposition}

\proof
Let $d=\gcdp{p}$. By Proposition \ref{reduced}, $r(p)=r(p/d)$, $d\,fw(p/d)=fw(p)$, while by definition of $p/d, 
\min(p)=d\,\min(p/d)$, $\max(p)=d\,\max(p/d)$ and $|p\mkern1mu|=|p/d|$. It suffices therefore to prove
the result for $p\in\ofs$ with $\gcdp{p}=1$, and this we shall do by induction on $\max(p)$. If $\gcdp{p}=1$
and $\max(p)=1$, then $p=(1)$ and so $fw(p)=1=\min(p)=\max(p)=|p\mkern1mu|$, while $r(p)=0$, 
so $\min(p)+\max(p)-(|p\mkern1mu|-1-r(p))=2\ge 1=fw(p)$.

Suppose now that $\gcdp{p}=1$, $\max(p)>1$, and the result holds for all elements of $\ofs$ with greatest
common divisor 1 and smaller maximum entry. Since $|p\mkern1mu|=1$ would imply that $1=\gcdp{p}=max(p)>1$, it
follows that $n=|p\mkern1mu|>1$. Consider $\hat{p}$, the totally reduced form of $p$. Since $\min(\hat{p})=p_1$,
$\max(\hat{p})=p_j$ for some $j$ with $2\le j\le n$, we have $\max(\hat{p})\le p_n$. If we are able to prove
that $fw(\hat{p})\le \min(\hat{p})+\max(\hat{p})-(|\hat{p}|-1)$, then by Corollary \ref{fw of totally reduced form},
$fw(p)=fw(\hat{p})\le \min(\hat{p})+\max(\hat{p})-(|\hat{p}|-1)\le p_1+p_n-(|p\mkern1mu|-r(p)-1)$, as required. Thus we
may assume that $p$ is totally reduced, and we are to prove that $fw(p)\le p_1+p_n-(|p\mkern1mu|-1)$.
Since $p$ is totally reduced, it is in particular trim, and so $fw(p)=f(p)=p_1+f(R(p))$. If $R(p)$ is not trim,
then by Proposition \ref{lower bound for fw}, $f(p)=2p_1$, and since $p_n-p_1\ge |p\mkern1mu|-1$, we have $p_1\le
p_n-(|p\mkern1mu|-1)$ and so $fw(p)=f(p)=2p_1\le p_1+p_n-(|p\mkern1mu|-1)$, as required. Thus we may assume 
that $R(p)$ is trim, so
$fw(R(p))=f(R(p))$. Furthermore, by Proposition \ref{first reduction}, the fact that $p$ is totally reduced means that
$p$ is reduced and so in particular, no entry of $p$ is a multiple of $p_1$. Thus $p_1\ne p_j-p_1$ for every
$j$ with $1\le j\le n$, so $|R(p)|=|p\mkern1mu|$. Apply the induction hypothesis to $R(p)$ to obtain 
$$\begin{array}{rl}
 fw(p)&=p_1+f(R(p))=p_1+fw(R(p))\\
 &\le p_1+\min(R(p))+\max(R(p))-(|R(p)|-1-r(R(p)))\\
 &=p_1+\min(R(p))+\max(R(p))-(|p\mkern1mu|-1-r(R(p))).\end{array}
$$
It will suffice to prove that $\min(R(p))+\max(R(p)) +r(R(p))\le p_n$. Let us first treat the case when $R(p)$ is
totally reduced; that is, $r(R(p))=0$. There are three subcases to consider. If $p_1\le p_2-p_1$, then $\min(R(p))=p_1$
and $\max(R(p))=p_n-p_1$, so $\min(R(p))+\max(R(p))+r(R(p))=p_1+p_n-p_1+0=p_n$, as required. Suppose now that
$p_2-p_1<p_1\le p_n-p_1$, so that $\min(R(p))=p_2-p_1$ and $\max(R(p))=p_n-p_1$. Then $\min(R(p))+\max(R(p))+r(R(p))
=p_2-p_1+p_n-p_1+0=p_n+p_2-2p_1$. But from $p_2-p_1<p_1$, we have $p_2-2p_1<0$ and so $p_n+p_2-2p_1<p_n$. Finally,
suppose that $p_n-p_1<p_1$, so that $\min(R(p))=p_2-p_1$ and $\max(R(p))=p_1$, which implies that 
$\min(R(p))+\max(R(p))+r(R(p))=p_2-p_1+p_1=p_2\le p_n$, as required.

We now treat the case when $R(p)$ is not totally reduced, so that $r(R(p))>0$. Let $j$ be such that $p_j-p_1$ is redundant
in $R(p)$. Let $S$, respectively $S'$, denote the initial segment of $R(p)$ that consists of the entries up to but not 
including $p_j-p_1$, respectively the entries up and including $p_j-p_1$, so $p_j-p_1\ge fw(S)$ and $\gcdp{S}=\gcdp{S'}$.
Since $p$ is totally reduced and therefore reduced, $p_1=p_j-p_1$ is not possible. Consider first the possibility
that $p_1<p_j-p_1$. Then either 
\begin{align*}
 S&=(p_2-p_1,p_3-p_1,\ldots,p_1,\ldots,p_{j-1}-p_1),\ \text{and}\\
 S'&= (p_2-p_1,p_3-p_1,\ldots,p_1,\ldots,p_{j}-p_1),
 \end{align*}
 or else
 \begin{align*}
  S&=(p_1,p_2-p_1,p_3-p_1,\ldots,p_{j-1}-p_1),\ \text{and}\\
  S'&= (p_1,p_2-p_1,p_3-p_1,\ldots,p_{j}-p_1),
 \end{align*}
 depending on whether $p_2-p_1< p_1$ or $p_2-p_1<p_1$. In either case, we have
 $S=R(p)\rest{j-1}=R(p\rest{j-1})$ and $S'=R(p\rest{j})$. Thus $R(p\rest{j})$ is not trim, so
 by Proposition \ref{lower bound for fw}, $f(p\rest{j})=2p_1$. Now $p$ is totally reduced, 
 so $p\rest{j}$ is trim, and since $J>1$, Proposition \ref{lower bound for fw} implies that
 $f(p\rest{j})>\max(p\rest{j})=p_j$. Thus $2p_1>p_j$; equivalently, $p_j-p_1<p_1$, contradicting our assumption that
$p_j-p_1>p_1$. Thus $p_j-p_1>p_1$ cannot hold, and since $p$ is reduced, $2p_1\ne p_j$ and thus we must have 
$p_j-p_1<p_1$. Since $j=2$ would mean that
$p_2-p_1$ is redundant and thus not the minimum entry of $R(p)$, it must be that $j>2$ and
so we have established that if $j$ is any index such that $p_j-p_1$ is  redundant in $R(p)$, then $p_2<p_j<2p_1$. Thus
$$
0<r(R(p))\le |\set j\mid p_2<p_j<2p_1\endset|+1,
$$
where we add 1 to acknowledge that $p_1$ might be redundant in $R(p)$. Thus $0<r(R(p))\le (2p_1-p_2-1)+1=2p_1-p_2$,
and so $p_2<2p_1$, which means that $p_2-p_1=\min(R(p))$. We consider two cases according to whether 
 $p_1<p_n-p_1$ or $p_1>p_n-p_1$ (again, $p_1=p_n-p_1$ is not possible since $p$ is totally reduced).
 
 Case 1: $p_1<p_n-p_1$. Then $\max(R(p))=p_n-p_1$, and so $\min(R(p))+\max(R(p))+r(R(p))\le p_2-p_1+p_n-p_1+2p_1-p_2=
 p_n$.
 
 Case 2: $p_1>p_n-p_1$. Then $\max(R(p))=p_1$, and $p_n<2p_1$. Recall that $R(p)$ is trim and so $\max(R(p))=p_1$
 is not redundant in $R(p)$. Thus in this case, we have
 $$
 r(R(p))\le |\set j\mid p_2<p_j<2p_1\endset|=|\set 3,4,\ldots,n\endset|=n-2.
 $$
 Thus $\min(R(p))+\max(R(p))+r(R(p))\le p_2-p_1+p_1+n-2=p_2+n-2\le p_n$. This completes the proof of the inductive
 step, and so the result follows by induction.
\endproof 

\begin{corollary}\label{gen of fine wilf}
 Let $p\in\ofs$ be totally reduced. Then $fw(p)\le \min(p)+\max(p)-\gcdp{p}(|p\mkern1mu|-1)$.
\end{corollary}

\begin{proposition}\label{truncated sequence fw}
 Let $p\in\ofs$ be such that $|p\mkern1mu|>1$, $\gcdp{p\rest {|p\mkern1mu|-1}}=\gcdp{p}$, and $\max(p)<fw(p\rest{|p\mkern1mu|-1})$.
 Then $fw(p\rest{|p\mkern1mu|-1})\ge fw(p)$.
\end{proposition}

\proof 
By Proposition \ref{reduced}, it suffices to prove this for $p\in\ofs$ for which $\gcdp{p}=1$.
Note that by Proposition \ref{min(p)=gcd(p)}, if $\min(p)=1$, then $fw(p)=1=fw(p\rest{|p\mkern1mu|-1})$, and so the result
holds in this case. Thus we may assume that $\min(p)>1=\gcdp{p}$.

Thus we consider $p\in\ofs$ such that $|p\mkern1mu|>1$, $\gcdp{p\rest {|p\mkern1mu|-1}}=1<\min(p)$, and
$\max(p)<fw(p\rest{|p\mkern1mu|-1})$. Suppose, contrary to our claim, that $fw(p\rest{|p\mkern1mu|-1})< fw(p)$.
Let $k= fw(p)-1-fw(p\rest{|p\mkern1mu|-1})$. 
If $k=0$, then $G'_p=G(p,fw(p\rest{|p\mkern1mu|-1}))$, which
has $G(p\rest{|p\mkern1mu|-1},fw(p\rest{|p\mkern1mu|-1}))$ as a spanning subgraph. By Corollary \ref{component count for G(p,k)} ,
$G(p\rest{|p\mkern1mu|-1},fw(p\rest{|p\mkern1mu|-1}))$ is connected. Thus if $k=0$, then  $G'_p$ is connected.
However, by Proposition \ref{optimality} (since $\min(p)>1=\gcdp{p}$), $G'_p$ has more than $\gcdp{p}=1$
components, so we would have a contradiction. Thus $k>0$. Now by Proposition \ref{upper bound for fw(p)} and
the fact that $\max(p)<fw(p\rest{|p\mkern1mu|-1})$, we have
$$\begin{array}{rl}
  k&=fw(p)-1-fw(p\rest{|p\mkern1mu|-1})\\
  &\le \min(p)+\max(p)-|p\mkern1mu|+1+r(p)-1-\max(p)\\
  &=\min(p)-(|p\mkern1mu|-r(p))<\min(p)<\max(p)<fw(p\rest{|p\mkern1mu|-1}).\end{array}
$$  
But every element of $\set 1,2,\ldots,fw(p\rest{|p\mkern1mu|-1})\endset$ is connected to 1 by a path in
$G(p\rest{|p\mkern1mu|-1},fw(p\rest{|p\mkern1mu|-1}))$. Since $\max(p)<fw(p\rest{|p\mkern1mu|-1})$, it follows 
that every element of $\set 
1,\ldots,fw(p)-1\endset$
is connected to 1 by a path in $G(p\rest{|p\mkern1mu|-1},fw(p)-1)$ and thus in $G_p'$, which contradicts the
fact that $G_p'$ is not connected.
\endproof

Note that it is possible that under the conditions of the preceding result, we may actually have $fw(p\rest{|p\mkern1mu|-1})> fw(p)$.
The lexically first such example would be $p=(5,7,8)$, where $fw(p)=10$, while $fw(5,7)=11$.

It might be tempting to believe that $fw$ grows monotically with respect to the product order on sequences of a given length and
greatest common divisor 1, and, as the Fine-Wilf theorem tells us, this is indeed the case for sequences of length 2 . However, 
this observation does not hold even for sequences of length 3. For example, $fw(7,9,11)=15$, while $fw(7,9,13)=14$.

\section{Combinatorics on words}

Let $p\in\ofs$ with $\gcdp{p}=1$ and $|p\mkern1mu|>1$, and consider the tableau for the computation of $f(p)$.
Let $m$ be minimum subject to the requirement that $|p^{(m)}|=1$. Then $p^{(m)}=(1)$, and $p^{(m-1)}=(1,2)$.
For each $i$ with $0\le i\le m$, we shall call $p^{(i)}$ a jump if $f(p^{(i)})=2p^{(i)}_1$, and in the tableau
for the computation of $f(p)$, we shall prefix each jump with a plus sign ($+$). Furthermore, let $J(p)$
denote the number of jumps in the tableau for the calculation of $f(p)$. For example, $p=(6,10,13)$ has
tableau

\centerline{\begin{tabular}{ll}
 &6,10,13\\
+&4,6,7\\
+&2,3,4\\
+&1,2\\
 &1\end{tabular}}

\noindent and so $J(p)=3$. We observe that $p^{(m)}$ is never a jump, while $p^{(m-1)}$ is always a jump. For each $i=1,\ldots,m$, let
$G^i$ denote the graph $G(p^{(m-i)},f(p^{(m-i)})-1)$, so that $G^1$ is the null graph on a single vertex, and $G^m=G(p,f(p)-1)$. 
Now if $m>1$; equivalently, $p\ne (1)$, for each $i=2,\ldots,m$, let $\alpha_i\from G^{i-1}\to G^i$ denote $\alpha_{p^{(m-(i-1))},f(p^{(m-(i-1))})-1}$
(see Definition \ref{component map}). This makes sense as 
$$
G(R(p^{(m-i)}),f(p^{(m-i+1)})-1)=G(p^{(m-i+1)},f(p^{(m-i+1)})-1)=G^{i-1}
$$
and
$$
 G(p^{(m-i)},p^{(m-i)}_1+f(p^{(m-i+1)})-1)=G(p^{(m-i)},f(p^{(m-i)})-1)=G^i
$$
and $\alpha_{p^{(m-(i-1))},f(p^{(m-(i-1))})-1}$ is a function from $G(R(p^{(m-i)}),f(p^{(m-i+1)})-1)$
to $G(p^{(m-i)},p^{(m-i)}_1+f(p^{(m-i+1)})-1)$.
Thus, for a vertex $j$, $\alpha_i(j)=p^{(m-i)}_1+j$. By Proposition \ref{alpha injective on components}, for each $i$, 
$\alpha_i$ induces an injective map from the set of components of $G^{i-1}$ into the set of components of $G^i$,
and $G^0$ has a single component. Moreover, the image of $\alpha_i$ is the set 
$$
\set p^{(m-i)}_1+1,\ldots,p^{(m-i)}_1+f(p^{(m-i+1)})-1\endset=\set p^{(m-i)}_1+1,\ldots,f(p^{(m-i)})-1\endset.
$$
If $f(p^{(m-i)})>2p^{(m-i)}_1$, then for
each $k\in \set 1,2,\ldots,p^{(m-i)}_1\endset$, $k+p^{(m-i)}_1\le 2p^{(m-i)}_1\le f(p^{(m-i)})-1$, and $\set k,k+p^{(m-i)}_1\endset$
is an edge in $G^i$ joining $k$ to a vertex in the image of $\alpha_i$, so $G^i$ and $G^{i-1}$ have exactly the same
number of components. On the other hand, if $f(p^{(m-i)})=2p^{(m-i)}_1$, then $p^{(m-i)}_1$ has degree 0 in $G^i$,
so $\set p^{(m-i)}_1\endset$ is a component of $G^i$ that is not contained in the image of $\alpha_i$. For any
$k$ with $1\le k<p^{(m-i)}$, $\set k,k+p^{(m-i)}\endset$ is an edge in $G^i$ joining $k$ to a vertex in the image of $\alpha_i$,
and so the number of components of $G^i$ is one more than the number of components of $G^{i-1}$. This proves the following
result.

\begin{proposition}\label{jumps=components}
 Let $p\in\ofs$ with $|p\mkern1mu|>1$ and $\gcdp{p}=1$. Then the number of components of $G(p,f(p)-1)$ is equal to $J(p)$, 
 the number of jumps in the tableau for the computation of $f(p)$.
\end{proposition}

Now, for $p=(p_1,p_2,\ldots,p_n)$ trim, $G(p,f(p)-1)=G_p'$, and the word $w$ of length $f(p)-1$ formed by labelling the 
$k$ components of $G_p'$ using, say, the integers from
$0$ to $k-1$, then setting $w_i$ equal to the label of the component containing vertex $i$ of $G_p'$ has periods $p_1,p_2,\ldots,p_n$,
but not $\gcdp{p}$. By Proposition \ref{jumps=components}, the number of distinct letters in the word is equal to $J(p)$.
Moreover, by Proposition \ref{palindrome}, $w$ is a palindrome.

We observe that the preceding discussion also shows that $w$ can be calculated from the tableau for the calculation of $f(p)$.
We begin at row $p^{(m-1)}$ with word $0$. Then at stage $p^{(i)}$, shift the preceding word $p^{(i)}_1$ spaces to the right.
If $p^{(i)}$ is not a jump, then the preceding word has length at least $p^{(i)}_1$ and we fill in the first $p^{(i)}_1$ locations of the new word with the first $p^{(i)}_1$ entries
in the preceding word, while if $p^{(i)}$ is a jump, then the preceding word has length $p^{(i)}_1-1$, and we fill in the
first $p^{(i)}_1-1$ spaces with the first $p^{(i)}_1-1$ entries of the preceding word and 
then introduce a new symbol for the vertex at position $p{(i)}_1$.

For example, $p=(6,10,13)$ has tableau

\centerline{\begin{tabular}{lll}
 &6,10,13 &0102010102010\\
+&4,6,7 &0102010\\
+&2,3,4 &010\\
+&1,2 &0\\
 &1&\end{tabular}}

\noindent and we have shown the construction of the word $w=010\,2\,01010\,2\,010$.
Note that $J(6,10,13)=3$ and indeed, $w$ has three distinct letters.

\end{document}